\documentclass[a4paper,oneside,11pt, reqno]{amsart}

\usepackage{mathrsfs}
\usepackage{amsfonts,amssymb,amsmath,amsthm}
\usepackage{url}
\usepackage{enumerate}
\usepackage{enumitem}
\usepackage{tikz}
\usepackage{float}
\usepackage{amscd}
\usepackage{amsxtra}
\usepackage[T1]{fontenc}
\usepackage[utf8]{inputenc}

\usetikzlibrary{arrows}
\usetikzlibrary{calc}
\usepackage{hyperref}

\theoremstyle{plain}
\newtheorem{theorem}{Theorem}[section]

\newtheorem{lemma}[theorem]{Lemma}
\newtheorem{proposition}[theorem]{Proposition}
\newtheorem{corollary}[theorem]{Corollary}

\theoremstyle{definition}
\newtheorem{definition}[theorem]{Definition}

\theoremstyle{remark}
\newtheorem{remark}[theorem]{Remark}

\newtheorem{case-new}{Case}

\renewcommand{\theenumi}{\roman{enumi}}
\renewcommand{\labelenumi}{(\theenumi )}
\numberwithin{equation}{section}

\usepackage{etoolbox}

\makeatletter
\@namedef{subjclassname@2020}{%
  \textup{2020} Mathematics Subject Classification}
\makeatother

\usepackage{epsfig}
\usepackage{mathrsfs}

\newcommand{\ncom}{\newcommand}
\ncom{\ul}{\underline}
\ncom{\ol}{\overline}
\ncom{\bq}{\begin{equation}}
\ncom{\eq}{\end{equation}}
\ncom{\beqn}{\begin{eqnarray*}}
\ncom{\eeqn}{\end{eqnarray*}}
\ncom{\beq}{\begin{equation}}
\ncom{\eeq}{\end{equation}}
\ncom{\nno}{\nonumber}
\ncom{\rar}{\rightarrow}
\ncom{\Rar}{\Rightarrow}
\ncom{\noin}{\noindent}
\ncom{\bc}{\begin{center}}
\ncom{\ec}{\end{center}}
\ncom{\sz}{\scriptsize}
\ncom{\rf}{\ref}
\ncom{\sgm}{\sigma}
\ncom{\Sgm}{\Sigma}
\ncom{\dt}{\delta}
\ncom{\Dt}{\Delta}
\ncom{\lmd}{\lambda}
\ncom{\Lmd}{\Lambda}
\ncom{\eps}{\epsilon}
\ncom{\pcc}{\stackrel{P}{>}}
\ncom{\dist}{{\rm\,dist}}
\ncom{\sspan}{{\rm\,span}}
\ncom{\im}{{\rm Im\,}}
\ncom{\sgn}{{\rm sgn\,}}
\ncom{\ba}{\begin{array}}
\ncom{\ea}{\end{array}}
\ncom{\eop}{\hfill{{\rule{2.5mm}{2.5mm}}}}
\ncom{\eoe}{\hfill{{\rule{1.5mm}{1.5mm}}}}
\ncom{\eof}{\hfill{{\rule{1.5mm}{1.5mm}}}}
\ncom{\hone}{\mbox{\hspace{1em}}}
\ncom{\htwo}{\mbox{\hspace{2em}}}
\ncom{\hthree}{\mbox{\hspace{3em}}}
\ncom{\hfour}{\mbox{\hspace{4em}}}
\ncom{\hsev}{\mbox{\hspace{7em}}}
\ncom{\vone}{\vskip 2ex}
\ncom{\vtwo}{\vskip 4ex}
\ncom{\vonee}{\vskip 1.5ex}
\ncom{\vthree}{\vskip 6ex}
\ncom{\vfour}{\vspace*{8ex}}
\ncom{\norm}{\|\;\;\|}
\ncom{\integ}[4]{\int_{#1}^{#2}\,{#3}\,d{#4}}
\ncom{\inp}[2]{\langle{#1},\,{#2} \rangle}
\ncom{\Inp}[2]{\Langle{#1},\,{#2} \Langle}
\ncom{\vspan}[1]{{{\rm\,span}\#1 \}}}
\ncom{\dm}[1]{\displaystyle {#1}}



\begin{document}
\title[Cyclic polynomials]{Cyclic polynomials in mixed-weighted Dirichlet-type spaces}

\author[R. Nailwal]{Rajkamal Nailwal}
\address{Rajkamal Nailwal, Nalanda University, Rajgir, Bihar, India}
\email{raj1994nailwal@gmail.com, rajkamal.2026211001@nalandauniv.edu.in}

\author[A. Zalar]{Alja\v z Zalar}
	\address{Alja\v z Zalar, 
        University of Ljubljana, 
		Faculty of Computer and Information Science \& 
		Faculty of Mathematics and Physics, \&
		Institute of Mathematics, Physics and Mechanics, Ljubljana, Slovenia.}
\email{aljaz.zalar@fri.uni-lj.si}
\thanks{AZ was supported by the Slovenian Research Agency core funding No. P1-0288 and grants J1-50002, J1-60011, J1-70017.}

\begin{abstract}
For $\boldsymbol{\alpha}=(\alpha_1,\alpha_2,\alpha_3)\in\mathbb R^3$,
we investigate cyclic polynomials in the mixed-weighted Dirichlet-type space $\mathscr{D}_{\boldsymbol{\alpha}}$ of holomorphic functions
$f(z_1,z_2)=\sum_{m,n\geq0}a_{m,n}z_1^mz_2^n$ on the unit bidisk satisfying
$$
\sum_{m,n\geq0}|a_{m,n}|^2
(m+1)^{\alpha_1}(n+1)^{\alpha_2}
(m+n+1)^{\alpha_3}<\infty.
$$
We give a complete characterization of cyclic polynomials in terms of
their zero sets and the parameters $\alpha_1,\alpha_2,\alpha_3$.
Our results extend previously known characterizations for
anisotropic and total-degree Dirichlet-type spaces on the unit bidisk. 
\end{abstract}
\maketitle
\tableofcontents
\section{Introduction}

The study of cyclic vectors in holomorphic function spaces lies at the
intersection of operator theory, approximation theory, and complex function
theory. In the classical Hardy space $H^2(\mathbb D),$ by Beurling's theorem, a function $f$ is
cyclic if and only if it is outer. In the classical Dirichlet space, the corresponding problem is substantially
more challenging. The celebrated Brown--Shields conjecture asserts that a function
in the Dirichlet space is cyclic if and only if it is outer and its boundary
zero set has logarithmic capacity zero \cite{BS1984}. This conjecture highlights the central
role played by boundary zero sets in cyclicity questions. The situation becomes
even more involved in several variables. Although cyclic polynomials in the
Hardy space $H^2(\mathbb D^d)$ are precisely those having no zeros in
$\mathbb D^d$ \cite{NGN1970}, their cyclicity in Dirichlet-type spaces may depend crucially
on both the geometry of their zero sets on the distinguished boundary and the
parameters defining the space. For a detailed discussion in one and several variables, we refer the reader to \cite{FKMR2014,FKR2009,JR2006, R1969}.

On the unit bidisk, an important family of examples is given by the anisotropic Dirichlet-type spaces
\[
\mathfrak D_{(\beta_1,\beta_2)}
=
\Big\{
\sum_{m,n\geq0}a_{m,n}z_1^mz_2^n:
\sum_{m,n\geq0}
(m+1)^{\beta_1}(n+1)^{\beta_2}|a_{m,n}|^2<\infty
\Big\}.
\]
where $\beta_1,\beta_2\in\mathbb R$.
B\'en\'eteau et al.~\cite{BKKLSS2016} obtained a complete
characterization of cyclic polynomials in the isotropic
product-weight spaces $\mathfrak D_{(\beta,\beta)}$,
$\beta\in\mathbb R$. This characterization was subsequently extended to the anisotropic setting by Knese et al.~in \cite{KKRS2019}.  

\begin{theorem}[{\cite{KKRS2019}}]
\label{thm:anisotropic-cyclicity}
Let $\beta_1,\beta_2\in\mathbb R$, and let
$p\in\mathbb C[z_1,z_2]$ be irreducible, depend on both
variables, and have no zeros in $\mathbb D^2$.
The following characterization holds in
$\mathfrak D_{(\beta_1,\beta_2)}$:
\begin{enumerate}
\item
If $\beta_1+\beta_2\leq1$, then $p$ is cyclic.

\item
If $\beta_1+\beta_2>1$ and
$\min\{\beta_1,\beta_2\}\leq1$, then $p$ is cyclic
if and only if $\mathcal Z(p)\cap\mathbb T^2$ is finite.

\item
If $\min\{\beta_1,\beta_2\}>1$, then $p$ is cyclic
if and only if $\mathcal Z(p)\cap\mathbb T^2=\varnothing$.
\end{enumerate}
\end{theorem}

A different scale of Dirichlet-type spaces on the unit bidisk is obtained by replacing the product weight with a total-degree weight:
$$
\mathcal D_{\gamma}
=
\Big\{
\sum_{m,n\geq 0}
a_{m,n}z_1^m z_2^n
:
\sum_{m,n\geq 0}
(m+n+1)^{\gamma}
|a_{m,n}|^2
<
\infty
\Big\},
$$
These spaces have been considered in \cite{BCG2025, BZ26, BS2026}. A complete characterization of cyclic polynomials in these spaces was obtained in  \cite{NZ26, Z25}. 

\begin{theorem}\label{thm:total-degree-cyclicity}
Let $\gamma\in\mathbb R$, and let
$p\in\mathbb C[z_1,z_2]$ be irreducible and have no zeros
in $\mathbb D^2$.
The following characterization holds in $\mathcal D_\gamma$:
\begin{enumerate}
\item
If $\gamma\leq1$, then $p$ is cyclic.

\item
If $1<\gamma\leq2$, then $p$ is cyclic
if and only if $\mathcal Z(p)\cap\mathbb T^2$ is finite.

\item
If $\gamma>2$, then $p$ is cyclic
if and only if $\mathcal Z(p)\cap\mathbb T^2=\varnothing$.
\end{enumerate}
\end{theorem}
Related results on cyclicity and invariant subspaces in
one-variable Dirichlet spaces can be found in
\cite{BS1984,FEK2016}.
In several variables, relevant work includes polynomial
cyclicity in Hardy spaces on polydisks \cite{NGN1970},
outer functions in function algebras on the bidisk
\cite{H88}, and cyclic polynomials in Dirichlet-type spaces
on the unit ball of $\mathbb C^2$ \cite{KV2023}.

The characterizations in
Theorems~\ref{thm:anisotropic-cyclicity}
and~\ref{thm:total-degree-cyclicity}
both relate cyclicity to the size of the zero set on
$\mathbb T^2$, but the parameter thresholds differ.
Consider an irreducible polynomial depending on both variables
and having no zeros in $\mathbb D^2$.
In the anisotropic product-weight setting, a polynomial with
infinitely many zeros on $\mathbb T^2$ is cyclic precisely
when $\beta_1+\beta_2\leq1$, whereas a polynomial with a
nonempty finite zero set on $\mathbb T^2$ is cyclic precisely
when $\min\{\beta_1,\beta_2\}\leq1$.
In the total-degree setting, the corresponding conditions are
$\gamma\leq1$ and $\gamma\leq2$, respectively.
These results lead to the question of whether the two
characterizations can be obtained as special cases of a more
general theorem.
This motivated us to study mixed-weighted Dirichlet-type spaces on the unit bidisk.

\subsection{Mixed-weighted Dirichlet-type spaces}

Let
$
\boldsymbol{\alpha}
=
(\alpha_1,\alpha_2,\alpha_3)
\in\mathbb R^3.
$
For $m,n\geq 0$, define
$$
w_{\boldsymbol{\alpha}}(m,n)
:=
(m+1)^{\alpha_1}
(n+1)^{\alpha_2}
(m+n+1)^{\alpha_3}.
$$

The mixed-weighted Dirichlet-type space
$\mathscr{D}_{\boldsymbol{\alpha}}$
consists of all holomorphic functions
$$
f(z_1,z_2)
=
\sum_{m,n\geq 0}
a_{m,n}z_1^m z_2^n
$$
on $\mathbb D^2$ satisfying
$$
\|f\|_{\boldsymbol{\alpha}}^2
:=
\sum_{m,n\geq 0}
w_{\boldsymbol{\alpha}}(m,n)
|a_{m,n}|^2
<
\infty.
$$
This is a reproducing kernel Hilbert space, with inner product
$$
\langle f,g\rangle_{\boldsymbol{\alpha}}
=
\sum_{m,n\geq 0}
w_{\boldsymbol{\alpha}}(m,n)
a_{m,n}\overline{b_{m,n}},\qquad
g(z_1,z_2)
=
\sum_{m,n\geq 0}
b_{m,n}z_1^m z_2^n.
$$
For $\zeta=(\zeta_1,\zeta_2)\in\mathbb D^2$, the reproducing
kernel at $\zeta$ is
\begin{equation}
\label{def:rep-kernel}
K_\zeta(z_1,z_2)
=
\sum_{m,n\geq0}
\frac{\overline{\zeta_1}^{\,m}\overline{\zeta_2}^{\,n}}
{w_{\boldsymbol{\alpha}}(m,n)}
z_1^mz_2^n.
\end{equation}
It belongs to $\mathscr{D}_{\boldsymbol{\alpha}}$ and satisfies
\[
f(\zeta)
=
\langle f,K_\zeta\rangle_{\boldsymbol{\alpha}},
\qquad f\in\mathscr{D}_{\boldsymbol{\alpha}}.
\]
The normalized monomials
\[
\Big\{
\frac{z_1^mz_2^n}{\sqrt{w_{\boldsymbol{\alpha}}(m,n)}}:
m,n\geq0
\Big\}
\]
form an orthonormal basis for 
$\mathscr{D}_{\boldsymbol{\alpha}}$.
In particular, polynomials are dense in
$\mathscr{D}_{\boldsymbol{\alpha}}$.

The mixed-weighted scale contains both previously studied
families: 
$$\mathscr{D}_{(\alpha_1,\alpha_2,0)}
=\mathfrak D_{(\alpha_1,\alpha_2)}\qquad \text{and}\qquad
\mathscr{D}_{(0,0,\alpha_3)}=\mathcal D_{\alpha_3}.$$
Also for fixed $\boldsymbol{\alpha}\in\mathbb R^3$, consider
the family of weights
\[
\mathfrak w_\gamma(m,n)
=
(m+1)^{\gamma\alpha_1}
(n+1)^{\gamma\alpha_2}
(m+n+1)^{(1-\gamma)\alpha_3},
\qquad 0\leq\gamma\leq1.
\]
The corresponding spaces are
$\mathscr{D}_{(\gamma\alpha_1,\gamma\alpha_2,
(1-\gamma)\alpha_3)}$.
This family continuously interpolates between the
total-degree and product weights:
\[
\begin{aligned}
\mathfrak w_0(m,n)
&=(m+n+1)^{\alpha_3},\\
\mathfrak w_1(m,n)
&=(m+1)^{\alpha_1}(n+1)^{\alpha_2},
\end{aligned}
\qquad m,n\geq0.
\]
Our results describe
how cyclicity changes as the weight varies continuously from a
total-degree weight to a product weight.
Moreover, we give a complete characterization of cyclic
polynomials in the mixed-weighted Dirichlet-type spaces
$\mathscr{D}_{\boldsymbol{\alpha}}$, recovering the known
anisotropic product-weight and total-degree characterizations
as special cases.

Since
$$
\sup_{m,n\geq 0}
\frac{
w_{\boldsymbol{\alpha}}(m+1,n)
}{
w_{\boldsymbol{\alpha}}(m,n)
}
<
\infty
\qquad
\text{and}
\qquad
\sup_{m,n\geq 0}
\frac{
w_{\boldsymbol{\alpha}}(m,n+1)
}{
w_{\boldsymbol{\alpha}}(m,n)
}
<
\infty,
$$
the coordinate multiplication operators
\beqn
M_{z_1}f(z_1,z_2)
=
z_1f(z_1,z_2) \qquad
\text{and} \qquad
M_{z_2}f(z_1,z_2)
=
z_2f(z_1,z_2)
\eeqn
are bounded on
$\mathscr{D}_{\boldsymbol{\alpha}}$.
Consequently, every polynomial is a multiplier of
$\mathscr{D}_{\boldsymbol{\alpha}}$.

\begin{definition}
For
$f\in\mathscr{D}_{\boldsymbol{\alpha}}$,
let
$$
[f]_{\boldsymbol{\alpha}}
=
\overline{
\operatorname{span}
\left\{
z_1^m z_2^n f:
m,n\geq 0
\right\}
}^{\,\mathscr{D}_{\boldsymbol{\alpha}}}.
$$
The function $f$ is said to be \emph{cyclic} in
$\mathscr{D}_{\boldsymbol{\alpha}}$
if
$$
[f]_{\boldsymbol{\alpha}}
=
\mathscr{D}_{\boldsymbol{\alpha}}.
$$
\end{definition}

We will use the following elementary observations about cyclicity.
\begin{enumerate}[leftmargin=*]
\item
Since polynomials are dense and multiplication by polynomials
is bounded, a function $f\in\mathscr{D}_{\boldsymbol{\alpha}}$
is cyclic if and only if there exists a sequence
$\{q_j\}_{j\geq1}\subseteq\mathbb C[z_1,z_2]$ such that
\[
\|q_jf-1\|_{\boldsymbol{\alpha}}\longrightarrow0.
\]

\item
Since polynomials are bounded multipliers of
$\mathscr{D}_{\boldsymbol{\alpha}}$, a factorization
$p=p_1\cdots p_k$ into nonconstant polynomials satisfies
\[
p\text{ is cyclic}
\iff
p_j\text{ is cyclic for every }j=1,\ldots,k.
\]
Thus, to characterize cyclic polynomials, it suffices
to consider irreducible polynomials.

\item
Since point evaluations on $\mathbb D^2$ are bounded,
a polynomial vanishing at a point of $\mathbb D^2$
cannot be cyclic.
Indeed, every polynomial multiple vanishes at that point
and therefore cannot approximate the constant function $1$
in norm.
Hence a necessary condition for cyclicity of $p$ is
\[
\mathcal Z(p)\cap\mathbb D^2=\varnothing.
\]
\end{enumerate}
Throughout the paper,
$\mathbb Z_+=\mathbb N_0=\{0,1,2,\ldots\}$.
For two nonnegative sequences $(u_{m,n})$ and $(v_{m,n})$,
we write $u_{m,n}\asymp v_{m,n}$ if there exist constants
$c,C>0$ such that
\[
c\,v_{m,n}\leq u_{m,n}\leq C\,v_{m,n}
\]
for all indices in the stated range.
Similarly, $u_{m,n}\lesssim v_{m,n}$ means that
$u_{m,n}\leq C\,v_{m,n}$, whereas
$u_{m,n}\gtrsim v_{m,n}$ means that
$u_{m,n}\geq c\,v_{m,n}$.
These constants may depend on fixed parameters, such as
$\boldsymbol{\alpha}$, but not on $m,n$.
We use the same notation for comparisons of norms, with
constants independent of the function being considered.

A finite zero set is allowed to be empty unless explicitly
stated otherwise.

\subsection{Three effective parameters}

The mixed weight naturally gives rise to three effective exponents:
\begin{equation}\label{eq:effective-parameters}
A=\alpha_1+\alpha_3,\qquad
B=\alpha_2+\alpha_3,\qquad
S=\alpha_1+\alpha_2+\alpha_3.
\end{equation}

Their geometric significance becomes clear from the restrictions of the weight to different regions of the lattice
$\mathbb Z_+^2$.
Along the first coordinate axis,
$$
w_{\boldsymbol{\alpha}}(m,0)
=
(m+1)^A.
$$
Along the second coordinate axis,
$$
w_{\boldsymbol{\alpha}}(0,n)
=
(n+1)^B.
$$
Along the diagonal,
$$
w_{\boldsymbol{\alpha}}(k,k)
=
(k+1)^{\alpha_1+\alpha_2}
(2k+1)^{\alpha_3}
\asymp
(k+1)^S.
$$

Thus, $A$ and $B$ describe the growth or decay of the weight
along the coordinate axes, whereas $S$ describes its behavior
when $m$ and $n$ grow at comparable rates.

\subsection{Main Result} 
 The following is the main result of this paper which gives a characterization of cyclic irreducible polynomials in mix-weighted Dirichlet-type spaces.

\begin{theorem}
\label{thm:main}
Let $\boldsymbol{\alpha}=(\alpha_1, \alpha_2, \alpha_3)\in\mathbb R^3$, set
$$
w_{\boldsymbol{\alpha}}(m,n)
=
(m+1)^{\alpha_1}(n+1)^{\alpha_2}(m+n+1)^{\alpha_3}, \qquad m,n \in \mathbb Z_+,
$$
and let $\mathscr{D}_{\boldsymbol{\alpha}}$ denote the corresponding
Dirichlet-type space on the bidisk. Set
\[
A=\alpha_1+\alpha_3,\qquad
B=\alpha_2+\alpha_3,\qquad
S=\alpha_1+\alpha_2+\alpha_3.
\]
Suppose that $p\in\mathbb C[z_1,z_2]$ is irreducible, depends on both
variables, and has no zeros in $\mathbb D^2$. 

\begin{enumerate}
    \item If $S\leq 1$, then $p$ is cyclic.

    \item If
    $
    S>1
   \,\text{and}\,
    \min\left\{A,B,\frac{S}{2}\right\}\leq1,
$
    then $p$ is cyclic if and only if it has finitely many zeros
    on the two-torus $\mathbb T^2$.

    \item If
    $
    \min\left\{A,B,\frac{S}{2}\right\}>1,
$
    then $p$ is cyclic if and only if it has no zeros on
    $\mathbb T^2$.
\end{enumerate}
\end{theorem}

Our main contribution is a complete characterization of cyclic
polynomials in $\mathscr{D}_{\boldsymbol{\alpha}}$ for every
$\boldsymbol{\alpha}\in\mathbb R^3$, including parameter ranges
not covered by direct norm comparisons with previously studied
spaces.
Two steps in the proof of Theorem~\ref{thm:main} are central
to this extension:
\begin{itemize}[leftmargin=*]
\item
For nonempty finite toral zero sets, we establish a sharp
cyclicity criterion for the model polynomial $2-z_1-z_2$.
We adapt the generating-function argument in
\cite[Appendix~A]{KKRS2019} and combine it with weighted
summability estimates for polynomial values on
$\mathbb N_0^2$.
We then use rotations of this model polynomial in the
smooth-quotient construction to obtain the finite-zero
classification. 

\item
For irreducible polynomials $p$ depending on both variables,
having no zeros in $\mathbb D^2$, and having infinitely many
zeros on $\mathbb T^2$, we prove noncyclicity whenever $S>1$
by constructing a nonzero bounded linear functional supported
on a regular arc of $\mathcal Z(p)\cap\mathbb T^2$.
This functional approach differs from the proof of
\cite[Theorem~5]{KKRS2019} in the anisotropic setting,
which reduces noncyclicity to the model polynomial $1-z_1z_2$
through algebraic branch substitutions and the symmetrization
lemma \cite[Lemma~8]{KKRS2019}.
\end{itemize}

The proofs of the remaining parts of the characterization in
Theorem~\ref{thm:main} adapt methods from
\cite{KKRS2019} and our earlier work \cite{NZ26}.

\subsection{Organization of the paper}

In Section~\ref{Prereq}, we characterize convergence of
$
\sum_{m,n\in\mathbb N_0}
\frac{1}{w_{\boldsymbol{\alpha}}(m,n)}.
$
This determines when evaluation at points of $\mathbb T^2$
is bounded and provides a noncyclicity criterion for
polynomials vanishing there.
We also establish norm comparisons between
$\mathscr{D}_{\boldsymbol{\alpha}}$ and anisotropic
Dirichlet-type spaces and discuss their consequences
for cyclicity.

In Section~\ref{cmp}, we determine the cyclicity of the model
polynomials
$1-z_1$,
$1-z_1z_2$,
$2-z_1-z_2$,
together with their rotations and monomial generalizations.
These results will be used in the proof of
Theorem~\ref{thm:main}.

Section~\ref{S-1} proves Theorem~\ref{thm:main}(i),
which establishes cyclicity when $S\leq1$.
In Section~\ref{p-h-f}, we characterize cyclicity for
irreducible polynomials depending on both variables,
having no zeros in $\mathbb D^2$, and having a nonempty
finite zero set on $\mathbb T^2$.

In Section~\ref{i-b-z}, we consider polynomials satisfying
the same irreducibility and zero-free assumptions but having
infinitely many zeros on $\mathbb T^2$.
We first construct a real-analytic arc of zeros whose slope
in angular coordinates is strictly negative and bounded
away from zero.
Integration by parts and Parseval's identity then give
a bounded linear functional supported on this arc,
proving noncyclicity when $S>1$.

In Section~\ref{sec:proof-main}, we combine the preceding
results to complete the proof of Theorem~\ref{thm:main}.

Finally, Section~\ref{sec:final-remarks} gives shorter
proofs of noncyclicity in two parameter ranges by comparison
with anisotropic Dirichlet-type spaces.
It also establishes a Hilbert-sum representation of
$\mathscr{D}_{\boldsymbol{\alpha}}$ when $\alpha_3\leq0$.

\section{Preliminaries} \label{Prereq}
In this section, we present some necessary conditions for cyclicity in $\mathscr{D}_{\boldsymbol{\alpha}}$. We also establish connections between $\mathscr{D}_{\boldsymbol{\alpha}}$ and anisotropic Dirichlet-type spaces.

\medskip

We call a polynomial \emph{stable} if it has no zeros in $\mathbb D^2$.
We first recall a geometric property of stable irreducible polynomials.

\begin{lemma}[{\cite[Lemma~3]{KKRS2019}}]
\label{lem:no-side-boundary-zeros}
Let $p\in\mathbb C[z_1,z_2]$ be irreducible, depend on both variables,
and have no zeros in $\mathbb D^2$. Then
$
\mathcal Z(p)
\cap
\bigl(
(\mathbb D\times\mathbb T)
\cup
(\mathbb T\times\mathbb D)
\bigr)
=
\varnothing.
$
Consequently, for every $r\in(0,1)$, the polynomial
$
p_r(z_1,z_2)
=
p(rz_1,z_2)
$
has no zeros on $\overline{\mathbb D}^{\,2}$.
\end{lemma}

\subsection{Bounded evaluation on the distinguished boundary}
For
$\tau=(\tau_1,\tau_2)\in\mathbb T^2$,
evaluation at $\tau$ extends continuously to
$\mathscr{D}_{\boldsymbol{\alpha}}$
if and only if
$$
\sum_{m,n\geq 0}
\frac{1}{
w_{\boldsymbol{\alpha}}(m,n)
}
<
\infty,
$$
since 
\beq \label{con-ser}
\left\|
K_\tau(z_1,z_2)
\right\|_{\boldsymbol{\alpha}}^2
=
\sum_{m,n\geq 0}
\frac{1}{
w_{\boldsymbol{\alpha}}(m,n)
},
\eeq
where 
$K_\tau(z_1,z_2)$ is as in \eqref{def:rep-kernel}.\\

 Next, we characterize convergence of the series in \eqref{con-ser}.

\begin{proposition} \label{finitesum}
Let
$\boldsymbol{\alpha}=(\alpha_1, \alpha_2, \alpha_3)\in\mathbb R^3$, and let $A,B,S$
be as in \eqref{eq:effective-parameters}.
Then
$$
\sum_{m,n\geq 0}
\frac{1}{
(m+1)^{\alpha_1}
(n+1)^{\alpha_2}
(m+n+1)^{\alpha_3}
}
<
\infty
$$
if and only if
$
A>1,
\,
B>1,$ and $
S>2.
$
\end{proposition}
\begin{proof}
We first split the sum into the regions
$
\Omega_1
=
\left\{
(m,n)\in\mathbb Z_+^2:
m> n
\right\}$
and $
\Omega_2
=
\left\{
(m,n)\in\mathbb Z_+^2:
n\geq m
\right\}.
$
On $\Omega_1$,
$$
m+n+1
\asymp
m+1.
$$
Therefore,
\begin{equation}\label{eq:omega1-series}
\sum_{(m,n)\in\Omega_1}
\frac{1}{w_{\boldsymbol{\alpha}}(m,n)}
\asymp
\sum_{m\geq0}(m+1)^{-A}
\sum_{n=0}^{m}(n+1)^{-\alpha_2}.
\end{equation}
The inner sum satisfies
\begin{equation}\label{eq:partial-sum-estimates}
\sum_{n=0}^{m}(n+1)^{-\alpha_2}
\asymp
\begin{cases}
1, & \alpha_2>1,\\[4pt]
\log(m+2), & \alpha_2=1,\\[4pt]
(m+1)^{1-\alpha_2}, & \alpha_2<1.
\end{cases}
\end{equation}
Indeed, comparison with an integral gives
\[
\sum_{n=0}^{m}(n+1)^{-\alpha_2}
\asymp
1+\int_1^{m+1}x^{-\alpha_2}\,dx.
\]
If $\alpha_2>1$, the integral is bounded independently of $m$,
so the right-hand side is comparable to $1$.
If $\alpha_2=1$, it equals $1+\log(m+1)$, which is comparable
to $\log(m+2)$.
Finally, if $\alpha_2<1$, then
\[
1+\int_1^{m+1}x^{-\alpha_2}\,dx
=
1+\frac{(m+1)^{1-\alpha_2}-1}{1-\alpha_2}
\asymp(m+1)^{1-\alpha_2}.
\]
In each case, the comparison constants are independent of $m$.

By \eqref{eq:partial-sum-estimates}, the series
\eqref{eq:omega1-series} converges if and only if $A>1$
when $\alpha_2\geq1$.
When $\alpha_2<1$, it is comparable to
\[
\sum_{m\geq0}(m+1)^{1-A-\alpha_2}
=
\sum_{m\geq0}(m+1)^{1-S},
\]
since $A+\alpha_2=S$.
Thus, in this case, \eqref{eq:omega1-series} converges
if and only if $S>2$.

By symmetry, the subseries over $\Omega_2$ converges
if and only if
\[
\begin{cases}
B>1, & \text{if }\alpha_1\geq1,\\[4pt]
S>2, & \text{if }\alpha_1<1.
\end{cases}
\]

Since $S=A+\alpha_2$, the condition for convergence on
$\Omega_1$ is equivalent to requiring both $A>1$ and $S>2$.
Indeed, if $\alpha_2\geq1$, then $A>1$ implies $S>2$;
if $\alpha_2<1$, then $S>2$ implies $A=S-\alpha_2>1$.
Similarly, convergence on $\Omega_2$ is equivalent to
requiring both $B>1$ and $S>2$.

Since $\Omega_1\cup\Omega_2=\mathbb Z_+^2$, the full series
converges if and only if both subseries converge. Therefore,
\[
\sum_{m,n\geq0}\frac{1}{w_{\boldsymbol{\alpha}}(m,n)}<\infty
\iff A>1,\ B>1,\ S>2.
\]
This completes the proof.
\end{proof}

The following consequence proves the necessity in Theorem~\ref{thm:main}(iii).

\begin{corollary}\label{ABS}
Let $\boldsymbol{\alpha}\in\mathbb R^3$,
and let $A,B,S$ be as in \eqref{eq:effective-parameters}.
Suppose that $A>1$, $B>1$, and $S>2$.
If $p\in\mathbb C[z_1,z_2]$ satisfies
$\mathcal Z(p)\cap\mathbb T^2\neq\varnothing$,
then $p$ is noncyclic in $\mathscr{D}_{\boldsymbol{\alpha}}$.
\end{corollary}

\begin{proof}
Choose
$
\tau\in\mathcal Z(p)\cap\mathbb T^2.
$
By Proposition~\ref{finitesum}, evaluation at $\tau$ is bounded. Since
$
(qp)(\tau)=0
$
for every polynomial $q$, the functional
$
f\longmapsto f(\tau)
$
annihilates
$[p]_{\boldsymbol{\alpha}}$.
However,
$
1(\tau)=1.
$
Therefore,
$
1\notin[p]_{\boldsymbol{\alpha}},
$
and hence $p$ is not cyclic.
\end{proof}

\begin{remark}
  By Proposition \ref{finitesum}, evaluation on
$\mathbb T^2$
fails to be bounded precisely when
$
A\leq 1
\,\text{or} \, 
B\leq 1\, \text{or}\,
S\leq 2.
$
Thus bounded point evaluations on $\mathbb T^2$ provide no obstruction to cyclicity in this parameter region. Theorem~\ref{fzc} shows that this condition is sufficient for cyclicity when the polynomial is irreducible, depends on both variables, is zero-free in $\mathbb D^2$, and has a nonempty finite toral zero set.
\end{remark}

\subsection{Relations with anisotropic Dirichlet-type spaces}

We compare the norm of $\mathscr{D}_{\boldsymbol{\alpha}}$
with those of anisotropic Dirichlet-type spaces to obtain
cyclicity and noncyclicity criteria for polynomials.

\begin{lemma}\label{lem:mixed-anisotropic-comparison}
Let $\boldsymbol{\alpha}\in\mathbb R^3$, and let $A,B$
be as in \eqref{eq:effective-parameters}.
Then the following statements hold.
\begin{enumerate}
\item
If $\alpha_3\leq0$, then
\[
\|f\|_{\boldsymbol{\alpha}}
\leq \|f\|_{\mathfrak D_{(A,\alpha_2)}},
\qquad
f\in\mathfrak D_{(A,\alpha_2)},
\]
and
\[
\|f\|_{\boldsymbol{\alpha}}
\leq \|f\|_{\mathfrak D_{(\alpha_1,B)}},
\qquad
f\in\mathfrak D_{(\alpha_1,B)}.
\]
In particular, $\mathfrak D_{(A,\alpha_2)}$ and
$\mathfrak D_{(\alpha_1,B)}$ embed continuously into
$\mathscr{D}_{\boldsymbol{\alpha}}$.

\item
If $\alpha_3>0$, then
\[
\mathscr{D}_{\boldsymbol{\alpha}}
=
\mathfrak D_{(A,\alpha_2)}
\cap
\mathfrak D_{(\alpha_1,B)}
\]
with equivalent norms:
\[
\|f\|_{\boldsymbol{\alpha}}^2
\asymp
\|f\|_{\mathfrak D_{(A,\alpha_2)}}^2
+
\|f\|_{\mathfrak D_{(\alpha_1,B)}}^2.
\]
\end{enumerate}
\end{lemma}

\begin{proof}
Suppose first that $\alpha_3\leq0$.
Since
\[
(m+n+1)^{\alpha_3}\leq(m+1)^{\alpha_3},
\qquad
(m+n+1)^{\alpha_3}\leq(n+1)^{\alpha_3},
\]
we have
\[
w_{\boldsymbol{\alpha}}(m,n)
\leq(m+1)^A(n+1)^{\alpha_2},
\qquad
w_{\boldsymbol{\alpha}}(m,n)
\leq(m+1)^{\alpha_1}(n+1)^B.
\]
Summing these inequalities against the squared Taylor
coefficients proves (i).

Suppose now that $\alpha_3>0$.
Put $M=\max\{m+1,n+1\}$.
The inequalities
\[
M\leq m+n+1\leq2M,
\qquad
M^{\alpha_3}
\leq(m+1)^{\alpha_3}+(n+1)^{\alpha_3}
\leq2M^{\alpha_3}
\]
give
\[
(m+n+1)^{\alpha_3}
\asymp
(m+1)^{\alpha_3}+(n+1)^{\alpha_3}.
\]
Multiplying by $(m+1)^{\alpha_1}(n+1)^{\alpha_2}$, we obtain
\[
w_{\boldsymbol{\alpha}}(m,n)
\asymp
(m+1)^A(n+1)^{\alpha_2}
+
(m+1)^{\alpha_1}(n+1)^B.
\]
Summing against the squared Taylor coefficients proves
the norm equivalence in (ii).
The mixed-weight norm is therefore finite if and only if
both anisotropic norms are finite, which also proves
the equality of spaces.
\end{proof}

\begin{remark}
Let $\boldsymbol{\alpha}\in\mathbb R^3$, and let $A,B,S$
be as in \eqref{eq:effective-parameters}.
Let $p\in\mathbb C[z_1,z_2]$ be irreducible, depend on both
variables, and have no zeros in $\mathbb D^2$.
The anisotropic spaces $\mathfrak D_{(A,\alpha_2)}$ and
$\mathfrak D_{(\alpha_1,B)}$ have the same parameter sum,
$A+\alpha_2=\alpha_1+B=S$.

\begin{enumerate}[leftmargin=*]
\renewcommand{\theenumi}{\roman{enumi}}
\renewcommand{\labelenumi}{\textup{(\theenumi)}}

\item
Suppose that $\alpha_3>0$, $S>1$, and $p$ has infinitely
many zeros on $\mathbb T^2$.
By Theorem~\ref{thm:anisotropic-cyclicity}(ii)--(iii),
$p$ is noncyclic in $\mathfrak D_{(A,\alpha_2)}$ and
$\mathfrak D_{(\alpha_1,B)}$.
If $p$ were cyclic in $\mathscr{D}_{\boldsymbol{\alpha}}$,
the norm comparison in
Lemma~\ref{lem:mixed-anisotropic-comparison}(ii)
would imply that $p$ is cyclic in
$\mathfrak D_{(A,\alpha_2)}$, a contradiction.
Thus $p$ is noncyclic in $\mathscr{D}_{\boldsymbol{\alpha}}$.

\item
Suppose that $\alpha_3>0$ and $S\leq1$.
By Theorem~\ref{thm:anisotropic-cyclicity}(i),
$p$ is cyclic in $\mathfrak D_{(A,\alpha_2)}$ and
$\mathfrak D_{(\alpha_1,B)}$.
This alone does not establish cyclicity in their intersection,
since the approximating sequences may differ.
By Lemma~\ref{lem:mixed-anisotropic-comparison}(ii),
cyclicity in $\mathscr{D}_{\boldsymbol{\alpha}}$ is equivalent
to the existence of a single sequence of polynomials $(q_j)$
such that
\[
\|pq_j-1\|_{\mathfrak D_{(A,\alpha_2)}}\longrightarrow0,
\qquad
\|pq_j-1\|_{\mathfrak D_{(\alpha_1,B)}}\longrightarrow0.
\]

\item
Suppose that $\alpha_3\leq0$ and $S\leq1$.
Since $A+\alpha_2=S$,
Theorem~\ref{thm:anisotropic-cyclicity}(i)
implies that $p$ is cyclic in $\mathfrak D_{(A,\alpha_2)}$.
By Lemma~\ref{lem:mixed-anisotropic-comparison}(i),
any sequence of polynomial multiples of $p$ converging to $1$
in $\mathfrak D_{(A,\alpha_2)}$ also converges to $1$ in
$\mathscr{D}_{\boldsymbol{\alpha}}$.
Hence $p$ is cyclic in $\mathscr{D}_{\boldsymbol{\alpha}}$.

\item
For $\alpha_3\leq0$, the inclusions in
Lemma~\ref{lem:mixed-anisotropic-comparison}(i)
can be strengthened to the identity
\[
\mathscr{D}_{\boldsymbol{\alpha}}
=
\mathfrak D_{(A,\alpha_2)}
+
\mathfrak D_{(\alpha_1,B)},
\]
with equivalence of Hilbert-sum norms.
This identification is established in
Proposition~\ref{prop:hilbert-sum}.
\end{enumerate}
\end{remark}

\section{Cyclicity of model polynomials}\label{cmp}

In this section, we study the cyclicity of model polynomials.
The results will be used in subsequent sections to prove parts of
Theorem~\ref{thm:main}.
The main result of this section is the following theorem.

\begin{theorem}\label{thm:model-polynomial-cyclicity}
Let $\boldsymbol{\alpha}=(\alpha_1,\alpha_2,\alpha_3)\in\mathbb R^3$,
$\zeta,\eta\in\mathbb T$, and let $j,k$ be positive integers.
Then the following statements hold.
\begin{enumerate}
\renewcommand{\theenumi}{\roman{enumi}}
\renewcommand{\labelenumi}{\textup{(\theenumi)}}
\item
$\zeta-z_1$ is cyclic in $\mathscr{D}_{\boldsymbol{\alpha}}$
if and only if $\alpha_1+\alpha_3\leq1$.

\item
$\eta-z_2$ is cyclic in $\mathscr{D}_{\boldsymbol{\alpha}}$
if and only if $\alpha_2+\alpha_3\leq1$.

\item
$1-z_1^jz_2^k$ is cyclic in $\mathscr{D}_{\boldsymbol{\alpha}}$
if and only if $\alpha_1+\alpha_2+\alpha_3\leq1$.

\item
$2-\overline{\zeta}z_1-\overline{\eta}z_2$ is cyclic in
$\mathscr{D}_{\boldsymbol{\alpha}}$ if and only if at least one
of the inequalities
\begin{equation}\label{eq:model-criterion}
\alpha_1+\alpha_3\leq1,\qquad
\alpha_2+\alpha_3\leq1,\qquad
\alpha_1+\alpha_2+\alpha_3\leq2
\end{equation}
holds.
\end{enumerate}
\end{theorem}

Taking $j=k=1$ in part~(iii) gives the criterion for $1-z_1z_2$.
We prove parts~(i)--(iii) by reduction to one variable.
For part~(iv), we study the coefficients of functions orthogonal
to polynomial multiples of $2-z_1-z_2$.
The auxiliary results below show that the corresponding arrays
have the form $b_{m,n}=P(m,n)$ for a polynomial
$P\in\mathbb C[x,y]$, and that the weighted summability
condition forces $P=0$ in the stated parameter region.

\subsection{Reduction to one variable}

For $\beta\in\mathbb R$, define the one-variable
Dirichlet-type space
\[
D_\beta
=
\Big\{
f=\sum_{m\geq0}c_mz^m\colon 
\|f\|_{D_\beta}^2
:=\sum_{m\geq0}(m+1)^\beta|c_m|^2<\infty
\Big\}.
\]
\begin{proof}[Proof of
Theorem~\ref{thm:model-polynomial-cyclicity}(i)--(iii)]
We first recall the one-variable criterion due to Brown and Shields \cite{BS1984}:
\begin{equation}\label{eq:one-variable-model-criterion}
1-z\text{ is cyclic in }D_\beta
\quad\Longleftrightarrow\quad \beta\leq1.
\end{equation}

For part~(i), rotation of the first variable reduces the proof
to $\zeta=1$. Given $q\in\mathbb C[z_1,z_2]$, write
\[
q(z_1,z_2)=q_0(z_1)+z_2q_1(z_1,z_2),
\]
where $q_0\in\mathbb C[z_1]$ and $q_1\in\mathbb C[z_1,z_2]$.
Orthogonality of monomials with distinct powers of $z_2$ gives
\[
\|(1-z_1)q-1\|_{\boldsymbol{\alpha}}^2
=\|(1-z_1)q_0-1\|_{\boldsymbol{\alpha}}^2
+\|(1-z_1)z_2q_1\|_{\boldsymbol{\alpha}}^2.
\]
Consequently,
\[
\inf_{q\in\mathbb C[z_1,z_2]}
\|(1-z_1)q-1\|_{\boldsymbol{\alpha}}
=
\inf_{q_0\in\mathbb C[z_1]}
\|(1-z_1)q_0-1\|_{\boldsymbol{\alpha}}.
\]
For $g(z_1)=\sum_{m\geq0}c_mz_1^m$, we have
\[
\|g\|_{\boldsymbol{\alpha}}^2
=\sum_{m\geq0}(m+1)^{\alpha_1+\alpha_3}|c_m|^2.
\]
It therefore suffices to determine whether $1-z$ is cyclic
in the one-variable space $D_{\alpha_1+\alpha_3}$.
Part~(i) now follows from
\eqref{eq:one-variable-model-criterion} with
$\beta=\alpha_1+\alpha_3$.
Interchanging the variables proves part~(ii).

For part~(iii), let $t=z_1^jz_2^k$.
Every polynomial $q$ has a unique decomposition
\[
q=q_0(t)+q_\perp,
\]
where $q_0\in\mathbb C[z]$ and $q_\perp$ contains no monomial
of the form $z_1^{j\ell}z_2^{k\ell}$, $\ell\geq0$.
Multiplication by $t$ preserves this property of $q_\perp$. Therefore,
\[
\|(1-t)q-1\|_{\boldsymbol{\alpha}}^2
=\|(1-t)q_0(t)-1\|_{\boldsymbol{\alpha}}^2
+\|(1-t)q_\perp\|_{\boldsymbol{\alpha}}^2.
\]
It follows that approximation of $1$ by polynomial multiples
of $1-t$ reduces to polynomials in $t$.
For $f(z)=\sum_{\ell\geq0}c_\ell z^\ell$, we have
\[
\|f(t)\|_{\boldsymbol{\alpha}}^2
=\sum_{\ell\geq0} (j\ell+1)^{\alpha_1}(k\ell+1)^{\alpha_2}
((j+k)\ell+1)^{\alpha_3}|c_\ell|^2.
\]
Since $j,k\geq1$ are fixed,
\[
(j\ell+1)^{\alpha_1}(k\ell+1)^{\alpha_2} ((j+k)\ell+1)^{\alpha_3}
\asymp(\ell+1)^{\alpha_1+\alpha_2+\alpha_3},
\]
with constants independent of $\ell$.
Hence cyclicity of $1-t$ in $\mathscr{D}_{\boldsymbol{\alpha}}$
is equivalent to cyclicity of $1-z$ in $D_{\alpha_1+\alpha_2+\alpha_3}$.
Part~(iii) now follows from
\eqref{eq:one-variable-model-criterion}.
\end{proof}

\begin{remark}
The zero set of $1-z_1z_2$ on $\mathbb T^2$ is
$\{(\zeta,\overline{\zeta}):\zeta\in\mathbb T\}$,
a curve on which both coordinates vary.
Part~(iii) shows that this polynomial is cyclic in
$\mathscr{D}_{\boldsymbol{\alpha}}$ precisely when
$\alpha_1+\alpha_2+\alpha_3\leq1$.
\end{remark}

\subsection{Polynomial solutions of the recurrence}

We now prepare the proof of
Theorem~\ref{thm:model-polynomial-cyclicity}(iv). Define
\[
p_0(z_1,z_2)=2-z_1-z_2.
\]
Its toral zero set consists of the single point $(1,1)$.
Our argument adapts the generating-function method of
\cite[Appendix~A]{KKRS2019}, also used in
\cite[Section~3]{NZ26}, to the mixed-weight setting.
We first express orthogonality to $[p_0]_{\boldsymbol{\alpha}}$
as a recurrence relation.

\begin{lemma}\label{lem:model-recurrence}
Let
\[
f(z_1,z_2) =
\sum_{m,n\geq 0}
\frac{b_{m,n}}{w_{\boldsymbol{\alpha}}(m,n)}z_1^mz_2^n
\in
\mathscr{D}_{\boldsymbol{\alpha}}.
\]
Then $ f\in[p_0]_{\boldsymbol{\alpha}}^\perp $ if and only if
\begin{equation} \label{eq:model-recurrence}
2b_{m,n} = b_{m+1,n}+b_{m,n+1},
\qquad
m,n\geq 0.
\end{equation}
Moreover,
\begin{equation} \label{eq:model-summability}
\sum_{m,n\geq 0}
\frac{|b_{m,n}|^2}{w_{\boldsymbol{\alpha}}(m,n)}
<
\infty.
\end{equation}
\end{lemma}

\begin{proof}
The summability condition \eqref{eq:model-summability} follows from
$f\in\mathscr{D}_{\boldsymbol{\alpha}}$:
\[
\|f\|_{\boldsymbol{\alpha}}^2
=
\sum_{m,n\geq 0}
w_{\boldsymbol{\alpha}}(m,n)
\left|
\frac{b_{m,n}}{w_{\boldsymbol{\alpha}}(m,n)}
\right|^2
=
\sum_{m,n\geq 0}
\frac{|b_{m,n}|^2}{w_{\boldsymbol{\alpha}}(m,n)}.
\]

By
\[
z_1^mz_2^np_0 = 2z_1^mz_2^n - z_1^{m+1}z_2^n - z_1^mz_2^{n+1},
\]
we obtain
\[
\left\langle
f,z_1^mz_2^np_0
\right\rangle_{\boldsymbol{\alpha}}
= 2b_{m,n} - b_{m+1,n} - b_{m,n+1}.
\]
Thus $ f\perp[p_0]_{\boldsymbol{\alpha}} $ if and only if
\eqref{eq:model-recurrence} holds for all
$m,n\geq 0$.
\end{proof}

The next lemma shows that every array satisfying the recurrence
\eqref{eq:model-recurrence} and the weighted summability condition
\eqref{eq:model-summability} consists of the values of a polynomial
on the nonnegative integer lattice.

\begin{lemma}
\label{lem:polynomial-recurrence}
Suppose that a complex array $\{b_{m,n}\}_{m,n\geq 0}$ satisfies
\begin{equation} \label{eq:polynomial-recurrence}
2b_{m,n} = b_{m+1,n}+b_{m,n+1},
\qquad m,n\geq 0,
\end{equation}
and
\begin{equation} \label{eq:polynomial-recurrence-summability}
\sum_{m,n\geq 0}
\frac{|b_{m,n}|^2}{w_{\boldsymbol{\alpha}}(m,n)}
<
\infty.
\end{equation}
Then there exists a polynomial $P\in\mathbb C[x,y]$ such that
\[
b_{m,n}=P(m,n),
\qquad m,n\geq 0.
\]
\end{lemma}

\begin{proof}
Define the generating function
\[
g(z_1,z_2)=\sum_{m,n\geq0}b_{m,n}z_1^mz_2^n.
\]
We divide the proof into three claims.

\medskip\noindent
\textbf{Claim 1.}
The series defining $g$ represents a holomorphic function on
$\mathbb D^2$, and there exists a polynomial $h$ such that
\[
g(z,0)=\frac{1}{z-1}
h\left(\frac{z}{z-1}\right),
\qquad z\in\mathbb D.
\]

\smallskip

\noindent\textit{Proof of Claim 1.}
We give the generating-function argument following
\cite[Appendix~A, pp.~44--45]{KKRS2019}.

Set
\[
\gamma=\sum_{j=1}^{3}\max\{\alpha_j,0\},
\qquad M=\max\{1,\gamma\}.
\]
Since $m+n+1\leq(m+1)(n+1)$, we have
\[
w_{\boldsymbol{\alpha}}(m,n)
\leq(m+n+1)^\gamma
\leq(m+1)^M(n+1)^M.
\]
Consequently, \eqref{eq:polynomial-recurrence-summability}
implies that
\[
K^2:=
\sum_{m,n\geq0}
\frac{|b_{m,n}|^2}{(m+1)^M(n+1)^M}
<\infty.
\]

Consider the generating function
\[
g(z_1,z_2)=\sum_{m,n\geq0}b_{m,n}z_1^mz_2^n.
\]
By the Cauchy--Schwarz inequality,
\[
\begin{aligned}
\sum_{m,n\geq0}|b_{m,n}|\,|z_1|^m|z_2|^n
&\leq
K\left(
\sum_{m,n\geq0}
(m+1)^M(n+1)^M|z_1|^{2m}|z_2|^{2n}
\right)^{1/2}\\
&=
K\left(\sum_{m\geq0}(m+1)^M|z_1|^{2m}\right)^{1/2}
 \left(\sum_{n\geq0}(n+1)^M|z_2|^{2n}\right)^{1/2}.
\end{aligned}
\]
The two series on the right converge uniformly on compact
subsets of $\mathbb D$. Hence the series defining $g$ converges
absolutely and locally uniformly on $\mathbb D^2$, and therefore
$g\in\operatorname{Hol}(\mathbb D^2)$.

Multiplying \eqref{eq:polynomial-recurrence} by
$z_1^{m+1}z_2^{n+1}$ and summing over $m,n\geq0$, we obtain
\[
2z_1z_2g(z_1,z_2)
=
z_2\bigl(g(z_1,z_2)-g(0,z_2)\bigr)
+
z_1\bigl(g(z_1,z_2)-g(z_1,0)\bigr).
\]
Thus
\[
(z_1+z_2-2z_1z_2)g(z_1,z_2)
=
z_1g(z_1,0)+z_2g(0,z_2).
\]

For $\zeta\neq\pm1$, put
\[
z_1=\frac{\zeta}{\zeta-1},
\qquad
z_2=\frac{\zeta}{\zeta+1}.
\]
A direct computation gives
\[
\left|\frac{\zeta}{\zeta-1}\right|<1
\iff \operatorname{Re}\zeta<\frac12,
\qquad
\left|\frac{\zeta}{\zeta+1}\right|<1
\iff \operatorname{Re}\zeta>-\frac12,
\]
and
\[
\frac{\zeta}{\zeta-1}
+
\frac{\zeta}{\zeta+1}
-
2\frac{\zeta}{\zeta-1}\frac{\zeta}{\zeta+1}
=0.
\]
Substituting into the generating-function identity yields
\[
\frac{\zeta}{\zeta-1}
g\left(\frac{\zeta}{\zeta-1},0\right)
+
\frac{\zeta}{\zeta+1}
g\left(0,\frac{\zeta}{\zeta+1}\right)
=0
\]
on the strip $-\frac12<\operatorname{Re}\zeta<\frac12$.
Dividing by $\zeta$ when $\zeta\neq0$, and then using
continuity at $\zeta=0$, gives
\[
\frac{1}{\zeta-1}
g\left(\frac{\zeta}{\zeta-1},0\right)
=
-\frac{1}{\zeta+1}
g\left(0,\frac{\zeta}{\zeta+1}\right).
\]
The two expressions are holomorphic on their respective
half-planes and agree on their intersection. Consequently,
\[
h(\zeta)=
\begin{cases}
\displaystyle
\frac{1}{\zeta-1}
g\left(\frac{\zeta}{\zeta-1},0\right),
& \operatorname{Re}\zeta<\frac12,\\[10pt]
\displaystyle
-\frac{1}{\zeta+1}
g\left(0,\frac{\zeta}{\zeta+1}\right),
& \operatorname{Re}\zeta>-\frac12
\end{cases}
\]
defines an entire function.

We next show that $h$ has polynomial growth.
Choose an integer $N\geq1$ such that $2N-1\geq M$.
For every $\ell\geq0$,
\[
(\ell+1)^M
\leq(\ell+1)^{2N-1}
\leq(2N-1)!\binom{\ell+2N-1}{2N-1}.
\]
Hence, for $0\leq r<1$,
\[
\sum_{\ell\geq0}(\ell+1)^M r^{2\ell}
\leq
\frac{(2N-1)!}{(1-r^2)^{2N}}.
\]
Restricting the definition of $K$ to the coordinate axes gives
\[
\sum_{m\geq0}\frac{|b_{m,0}|^2}{(m+1)^M}\leq K^2,
\qquad
\sum_{n\geq0}\frac{|b_{0,n}|^2}{(n+1)^M}\leq K^2.
\]
Thus, by the Cauchy--Schwarz inequality and the preceding estimate,
\[
\begin{aligned}
|g(z,0)|
&\leq
\left(\sum_{m\geq0}\frac{|b_{m,0}|^2}{(m+1)^M}\right)^{1/2}
\left(\sum_{m\geq0}(m+1)^M|z|^{2m}\right)^{1/2}\\
&\leq
\frac{K\sqrt{(2N-1)!}}{(1-|z|^2)^N}.
\end{aligned}
\]
The same estimate holds for $g(0,z)$. Consequently,
\[
|g(z,0)|+|g(0,z)|
\leq\frac{C}{(1-|z|^2)^N},
\qquad z\in\mathbb D,
\]
where $C$ is independent of $z$.

If $\operatorname{Re}\zeta\leq0$, then
\[
1-\left|\frac{\zeta}{\zeta-1}\right|^2
=
\frac{1-2\operatorname{Re}\zeta}{|\zeta-1|^2}
\geq\frac{1}{|\zeta-1|^2}.
\]
Using the first representation of $h$, we obtain
\[
\begin{aligned}
|h(\zeta)|
&\leq
\frac{C}{|\zeta-1|}
\left(
1-\left|\frac{\zeta}{\zeta-1}\right|^2
\right)^{-N}\\
&\leq C|\zeta-1|^{2N-1}
\leq C(1+|\zeta|)^{2N-1}.
\end{aligned}
\]
Similarly, if $\operatorname{Re}\zeta\geq0$, then
\[
1-\left|\frac{\zeta}{\zeta+1}\right|^2
=
\frac{1+2\operatorname{Re}\zeta}{|\zeta+1|^2}
\geq\frac{1}{|\zeta+1|^2}.
\]
The second representation of $h$ gives
\[
|h(\zeta)|
\leq C|\zeta+1|^{2N-1}
\leq C(1+|\zeta|)^{2N-1}.
\]
Thus
\[
|h(\zeta)|\leq C(1+|\zeta|)^{2N-1},
\qquad \zeta\in\mathbb C.
\]
For every integer $k>2N-1$, Cauchy's estimates imply
\[
|h^{(k)}(0)|
\leq
\frac{k!}{R^k}\max_{|\zeta|=R}|h(\zeta)|
\leq
Ck!\frac{(1+R)^{2N-1}}{R^k}
\longrightarrow0
\]
as $R\to\infty$. Therefore $h$ is a polynomial of degree
at most $2N-1$.

Finally, for $z\in\mathbb D$, the substitution
$\zeta=z/(z-1)$ in the first representation of $h$ gives
\[
g(z,0)
=
\frac{1}{z-1}
h\left(\frac{z}{z-1}\right).
\]
Likewise, substituting $\zeta=z/(1-z)$ into the second
representation gives
\[
g(0,z)
=
\frac{1}{z-1}
h\left(\frac{z}{1-z}\right).
\]
This proves Claim~1.
\hfill $\blacksquare$

\medskip

\noindent
\textbf{Claim 2.}
There exists a polynomial $P_0\in\mathbb C[x]$ such that
\[
b_{m,0}=P_0(m),
\qquad m\geq0.
\]

\smallskip

\noindent\textit{Proof of Claim 2.}
Write $h(\xi)=\sum_{j=0}^{d}c_j\xi^j$. By Claim~1,
\[
g(z,0)
=\frac{1}{z-1}h\left(\frac{z}{z-1}\right)
=\sum_{j=0}^{d}c_j\frac{z^j}{(z-1)^{j+1}}.
\]
For each $j\geq0$, the binomial-series identity gives
\[
\frac{1}{(1-z)^{j+1}}
=\sum_{\ell\geq0}\binom{\ell+j}{j}z^\ell,
\qquad |z|<1.
\]
Multiplying by $(-1)^{j+1}z^j$ and setting $m=\ell+j$,
we obtain
\[
\frac{z^j}{(z-1)^{j+1}}
=(-1)^{j+1}\sum_{m\geq j}\binom{m}{j}z^m.
\]
Substituting these expansions into the expression for $g(z,0)$
and collecting coefficients, we find
\[
g(z,0)
=\sum_{m\geq0}
\left(
\sum_{j=0}^{d}(-1)^{j+1}c_j\binom{m}{j}
\right)z^m,
\]
where $\binom{m}{j}=0$ for $m<j$.
Comparing this with $g(z,0)=\sum_{m\geq0}b_{m,0}z^m$ gives
\[
b_{m,0}
=\sum_{j=0}^{d}(-1)^{j+1}c_j\binom{m}{j}.
\]
Thus
\[
P_0(x)=\sum_{j=0}^{d}(-1)^{j+1}c_j\binom{x}{j}
\]
has the required values. Here
$\binom{x}{0}=1$,
$\binom{x}{j}=\frac{x(x-1)\cdots(x-j+1)}{j!}$,
$j\geq1,$
so $P_0$ is indeed a polynomial. This proves Claim~2.\hfill$\blacksquare$

\medskip\noindent
\textbf{Claim 3.}
There exists a polynomial $P\in\mathbb C[x,y]$ such that
\[
b_{m,n}=P(m,n),
\qquad m,n\geq0.
\]

\smallskip
\noindent\textit{Proof of Claim 3.}
If $P_0=0$, \eqref{eq:polynomial-recurrence} implies
inductively that $b_{m,n}=0$ for all $m,n\geq0$,
and we take $P=0$.
Otherwise, put $d_0=\deg P_0$ and define the forward-difference
operator $\Delta:\mathbb C[x]\to\mathbb C[x]$ by
\[
(\Delta Q)(x)=Q(x+1)-Q(x),
\qquad Q\in\mathbb C[x].
\]
Rewriting \eqref{eq:polynomial-recurrence} as
\[
b_{m,n+1}=2b_{m,n}-b_{m+1,n},
\]
we first prove that increasing the second index corresponds
to applying $I-\Delta$ to the polynomial in the first index.
Let
\begin{equation}\label{eq:forward-difference-identity}
((I-\Delta)Q)(x)=2Q(x)-Q(x+1),
\qquad Q\in\mathbb C[x].
\end{equation}
where $I$ denotes the identity operator on $\mathbb C[x]$.
We show by induction on $n$ that
\begin{equation}\label{eq:iterated-difference-representation}
b_{m,n}=\bigl((I-\Delta)^nP_0\bigr)(m),
\qquad m,n\geq0.
\end{equation}
For $n=0$, this is the identity established in Claim~2.
Suppose that the formula holds for some $n\geq0$ and every
$m\geq0$. By \eqref{eq:polynomial-recurrence},
\[
\begin{aligned}
b_{m,n+1}
&=2b_{m,n}-b_{m+1,n}\\
&=2\bigl((I-\Delta)^nP_0\bigr)(m)
  -\bigl((I-\Delta)^nP_0\bigr)(m+1)\\
&=\bigl((I-\Delta)((I-\Delta)^nP_0)\bigr)(m)\\
&=\bigl((I-\Delta)^{n+1}P_0\bigr)(m).
\end{aligned}
\]
In the third equality, we applied
\eqref{eq:forward-difference-identity} with
$Q=(I-\Delta)^nP_0$ and $x=m$. This completes the induction.

Each application of $\Delta$ lowers the degree of a
nonconstant polynomial by one and annihilates constants.
Hence $\Delta^jP_0=0$ for $j>d_0$.
Applying the binomial theorem to $(I-\Delta)^n$ in
\eqref{eq:iterated-difference-representation}, we obtain
\[
b_{m,n}
=\sum_{j=0}^{d_0}(-1)^j
\binom{n}{j}(\Delta^jP_0)(m)
\]
with $\binom{n}{j}=0$ when $j>n$.
Therefore,
\[
P(x,y)=\sum_{j=0}^{d_0}(-1)^j
\binom{y}{j}(\Delta^jP_0)(x)
\]
is a polynomial satisfying $P(m,n)=b_{m,n}$ for every
$m,n\geq0$. This proves Claim~3 and completes the proof of the lemma.
\end{proof}

\subsection{Divergence of weighted polynomial sums}

To prove the sufficiency direction of
Theorem~\ref{thm:model-polynomial-cyclicity}(iv), it suffices
to show that $[2-z_1-z_2]_{\boldsymbol{\alpha}}^\perp=\{0\}$.
By Lemmas~\ref{lem:model-recurrence}
and~\ref{lem:polynomial-recurrence}, this follows if, under
the stated parameter conditions, no nonzero polynomial $P$ satisfies
\[
\sum_{m,n\geq0}
\frac{|P(m,n)|^2}{w_{\boldsymbol{\alpha}}(m,n)}<\infty.
\]
The next lemma establishes this fact.

\begin{lemma}\label{lem:polynomial-array-divergence}
Let $\boldsymbol{\alpha}\in\mathbb R^3$ and let
$P\in\mathbb C[x,y]$ be nonzero. If at least one of the inequalities
\[
\alpha_1+\alpha_3\leq1,\qquad
\alpha_2+\alpha_3\leq1,\qquad
\alpha_1+\alpha_2+\alpha_3\leq2
\]
holds, then
\[
\sum_{m,n\geq0}
\frac{|P(m,n)|^2}{w_{\boldsymbol{\alpha}}(m,n)}=\infty.
\]
\end{lemma}

\begin{proof}
Suppose first that $\alpha_1+\alpha_3\leq1$.
Choose a nonnegative integer $n_0$ such that
$P(\cdot,n_0)$ is a nonzero polynomial, and let its degree be $d_0$.
Such an $n_0$ exists because a nonzero polynomial cannot vanish
identically on every horizontal line with nonnegative integer height.
For all sufficiently large $m$,
\[
|P(m,n_0)|\gtrsim m^{d_0},\qquad w_{\boldsymbol{\alpha}}(m,n_0)
\asymp m^{\alpha_1+\alpha_3}.
\]
Thus, for some sufficiently large integer $M$,
\[
\sum_{m\geq M}
\frac{|P(m,n_0)|^2}{w_{\boldsymbol{\alpha}}(m,n_0)}
\gtrsim\sum_{m\geq M}m^{2d_0-\alpha_1-\alpha_3}
=\infty,
\]
since $2d_0-\alpha_1-\alpha_3\geq-1$.
The case $\alpha_2+\alpha_3\leq1$ follows by symmetry.

Suppose next that $\alpha_1+\alpha_2+\alpha_3\leq2$.
Let $d=\deg P$ and let $P_d$ be its homogeneous part of degree $d$.
The polynomial $P_d(1,t)$ is nonzero, so there exist
$0<c_1<c_2$ and $\delta>0$ such that
\[
|P_d(1,t)|\geq\delta,\qquad t\in[c_1,c_2].
\]
Since $P-P_d$ has degree at most $d-1$ and
$n/m\in[c_1,c_2]$, there is a constant $C>0$, independent
of $m,n$, such that
\[
\left|m^{-d}P(m,n)-P_d(1,n/m)\right|
\leq \frac{C}{m},
\qquad m\geq1,\quad c_1m\leq n\leq c_2m.
\]
Thus the error tends to zero as $m\to\infty$, with the
same bound for every integer $n$ in the stated range.
Hence
$|P(m,n)|\gtrsim m^d$ on this region for sufficiently large $m$. Moreover,
\[
w_{\boldsymbol{\alpha}}(m,n)
\asymp m^{\alpha_1+\alpha_2+\alpha_3}
\quad\text{when }c_1m\leq n\leq c_2m.
\]
For all sufficiently large $m$,
\[
\#\{n\in\mathbb Z:c_1m\leq n\leq c_2m\}
=\lfloor c_2m\rfloor-\lceil c_1m\rceil+1
\geq \frac{c_2-c_1}{2}\,m.
\]
Consequently, choosing $M$ sufficiently large, we obtain
\[
\begin{aligned}
\sum_{m,n\geq0}
\frac{|P(m,n)|^2}{w_{\boldsymbol{\alpha}}(m,n)}
&\geq
\sum_{m\geq M}
\sum_{\lceil c_1m\rceil\leq n\leq\lfloor c_2m\rfloor}
\frac{|P(m,n)|^2}{w_{\boldsymbol{\alpha}}(m,n)}\\
&\gtrsim
\sum_{m\geq M}
\sum_{\lceil c_1m\rceil\leq n\leq\lfloor c_2m\rfloor}
m^{2d-\alpha_1-\alpha_2-\alpha_3}\\
&=
\sum_{m\geq M}
\bigl(\lfloor c_2m\rfloor-\lceil c_1m\rceil+1\bigr)
m^{2d-\alpha_1-\alpha_2-\alpha_3}\\
&\gtrsim
\sum_{m\geq M}
m^{2d+1-\alpha_1-\alpha_2-\alpha_3}
=\infty.
\end{aligned}
\]
This completes the proof of the lemma.
\end{proof}

\subsection{Completion of the proof of Theorem~\ref{thm:model-polynomial-cyclicity}}

\begin{proof}[Proof of
Theorem~\ref{thm:model-polynomial-cyclicity}(iv)]
Rotations of the variables preserve the norm of
$\mathscr{D}_{\boldsymbol{\alpha}}$ and send polynomial multiples
to polynomial multiples. It therefore suffices to consider $p_0=2-z_1-z_2$.

Suppose that none of the inequalities in
\eqref{eq:model-criterion} holds.
By Proposition~\ref{finitesum},
\[
\sum_{m,n\geq0}\frac{1}{w_{\boldsymbol{\alpha}}(m,n)}<\infty.
\]
Hence by the bounded-evaluation
criterion established before Proposition~\ref{finitesum}, evaluation at $(1,1)$
defines a bounded linear functional on
$\mathscr{D}_{\boldsymbol{\alpha}}$.
It vanishes on every polynomial multiple of $p_0$, but not
on the constant function $1$, so $p_0$ is noncyclic.

Conversely, suppose that at least one inequality in
\eqref{eq:model-criterion} holds.
Let $f\in[p_0]_{\boldsymbol{\alpha}}^\perp$ and write
\[
f(z_1,z_2)=\sum_{m,n\geq0}
\frac{b_{m,n}}{w_{\boldsymbol{\alpha}}(m,n)}z_1^mz_2^n.
\]
By Lemma~\ref{lem:model-recurrence}, the array $(b_{m,n})$
satisfies the recurrence \eqref{eq:model-recurrence} and
the summability condition \eqref{eq:model-summability}.
Lemma~\ref{lem:polynomial-recurrence} then gives
$b_{m,n}=P(m,n)$ for a polynomial $P\in\mathbb C[x,y]$.
By Lemma~\ref{lem:polynomial-array-divergence}, the summability
condition forces $P=0$.
Hence $f=0$, so $[p_0]_{\boldsymbol{\alpha}}^\perp=\{0\}$.
It follows that $p_0$ is cyclic, completing the proof.
\end{proof}

\section{The case \texorpdfstring{$S\leq 1$}{S <= 1}
of Theorem~\ref{thm:main}}\label{S-1}

In this section, we prove Theorem~\ref{thm:main}(i).
We begin with a standard criterion for weak convergence in
reproducing kernel Hilbert spaces.
For a Hilbert space $\mathcal H$, we write
$f_k\rightharpoonup f$ if
\[
\langle f_k,g\rangle_{\mathcal H}
\longrightarrow
\langle f,g\rangle_{\mathcal H}
\qquad\text{for every }g\in\mathcal H.
\]

\begin{lemma}\label{lem:bounded-pointwise-weak-convergence}
Let $\mathcal H$ be a reproducing kernel Hilbert space of functions on
a set $X$. Suppose that $\{f_k\}_{k\geq 1}\subseteq \mathcal H$ satisfies
$
\sup_{k\geq 1}\|f_k\|_{\mathcal H}<\infty
$
and
$
f_k(x)\longrightarrow f(x),
\,
x\in X,
$
for some $f\in \mathcal H$. Then
$
f_k\rightharpoonup f$ in $\mathcal H$.
\end{lemma}

\begin{proof}
For $x\in X$, let $K_x$ denote the reproducing kernel at $x$. Then
$$
\langle f_k,K_x\rangle_{\mathcal H}
=
f_k(x)
\longrightarrow
f(x)
=
\langle f,K_x\rangle_{\mathcal H}.
$$
Hence
$
\langle f_k,g\rangle_{\mathcal H}
\longrightarrow
\langle f,g\rangle_{\mathcal H}
$
for every
$
g\in
\operatorname{span}\{K_x:x\in X\}.
$
By \cite[Proposition~2.1]{PaulsenRaghupathi2016},
the linear span of $\{K_x:x\in X\}$ is dense in $\mathcal H$.
Since $\sup_{k\geq1}\|f_k\|_{\mathcal H}<\infty$,
the convergence extends to every $g\in\mathcal H$.
\end{proof}

We can now prove Theorem~\ref{thm:main}(i).

\begin{theorem}\label{s1}
Let $\boldsymbol{\alpha}=(\alpha_1,\alpha_2,\alpha_3)\in\mathbb R^3$,
and $A,B,S$ be as in \eqref{eq:effective-parameters}.
Let $p\in\mathbb C[z_1,z_2]$ be irreducible, depend on both
variables, and have no zeros in $\mathbb D^2$.
If $S\leq1$, then $p$ is cyclic in
$\mathscr{D}_{\boldsymbol{\alpha}}$.
\end{theorem}

\begin{proof}
Choose $r_k\uparrow1$ and set
\[
F_k(z_1,z_2)
=\frac{p(z_1,z_2)}{p(r_kz_1,z_2)}.
\]
Since $A+\alpha_2=\alpha_1+B=S\leq1$,
\cite[Lemma~2]{KKRS2019} gives $F_k\rightharpoonup1$
in both $\mathfrak D_{(A,\alpha_2)}$ and
$\mathfrak D_{(\alpha_1,B)}$.
Since weakly convergent sequences are bounded in norm
by the uniform boundedness principle, we obtain
\[
\sup_k\|F_k\|_{\mathfrak D_{(A,\alpha_2)}}<\infty,
\qquad
\sup_k\|F_k\|_{\mathfrak D_{(\alpha_1,B)}}<\infty.
\]

For every $\alpha_3\in\mathbb R$,
Lemma~\ref{lem:mixed-anisotropic-comparison}(i)--(ii)
gives
\[
\|F_k\|_{\boldsymbol{\alpha}}^2
\lesssim
\|F_k\|_{\mathfrak D_{(A,\alpha_2)}}^2
+
\|F_k\|_{\mathfrak D_{(\alpha_1,B)}}^2.
\]
Thus $(F_k)$ is bounded in
$\mathscr{D}_{\boldsymbol{\alpha}}$.
Since $F_k\to1$ pointwise on $\mathbb D^2$,
Lemma~\ref{lem:bounded-pointwise-weak-convergence}
implies that $F_k\rightharpoonup1$ in
$\mathscr{D}_{\boldsymbol{\alpha}}$.

We conclude as in the proof of
\cite[Theorem~4]{KKRS2019}. Namely, 
by Lemma~\ref{lem:no-side-boundary-zeros},
$1/p(r_kz_1,z_2)$ is holomorphic on a neighborhood of
$\overline{\mathbb D}^{\,2}$.
Its Taylor coefficients decay exponentially by Cauchy's
estimates, whereas $w_{\boldsymbol{\alpha}}(m,n)$ grows
at most polynomially in $m+n$.
Consequently, $1/p(r_kz_1,z_2)$ belongs to
$\mathscr{D}_{\boldsymbol{\alpha}}$.
Polynomial approximation of this reciprocal and boundedness
of multiplication by $p$ show that
$F_k\in[p]_{\boldsymbol{\alpha}}$.
Since $[p]_{\boldsymbol{\alpha}}$ is weakly closed,
its weak limit $1$ also belongs to
$[p]_{\boldsymbol{\alpha}}$.
Therefore, $p$ is cyclic.
\end{proof}

\section{Polynomials having finitely many zeros on \texorpdfstring{$\mathbb T^2$}{the two-torus}}\label{p-h-f}

In this section, we characterize the cyclicity of polynomials
in $\mathscr{D}_{\boldsymbol{\alpha}}$ having finitely many zeros on $\mathbb T^2$.
We first construct a quotient belonging to
$\mathscr{D}_{\boldsymbol{\alpha}}$.
\begin{lemma}\label{lem:smooth-quotient}
Let $p\in\mathbb C[z_1,z_2]$ be irreducible, depend on both
variables, and have no zeros in $\mathbb D^2$. Suppose that
\[
\mathcal Z(p)\cap\mathbb T^2
=
\{(\zeta_j,\eta_j):j=1,\ldots,N\},
\qquad N\geq1,
\]
and define
\[
L_j(z_1,z_2)
=
2-\overline{\zeta_j}z_1-\overline{\eta_j}z_2.
\]
Then, for every $\boldsymbol{\alpha}\in\mathbb R^3$, there exists
an integer $M\geq1$ such that
\[
Q_M
=
\frac{\prod_{j=1}^{N}L_j^M}{p}
\in\mathscr{D}_{\boldsymbol{\alpha}}.
\]
\end{lemma}

\begin{proof}
Choose an integer $s\geq1$ such that
\[
2s\geq
\max\{\alpha_1^++\alpha_3^+,\,
      \alpha_2^++\alpha_3^+\},
\qquad t^+=\max\{t,0\}.
\]
By \cite[Lemma~4.2]{NZ26}, for sufficiently large $M$,
the boundary quotient $Q_M$ extends to a function in
$C^{2s}(\mathbb T^2)$.
Then, the Fourier coefficients of $Q_M$ satisfy
\[
\sum_{m,n\geq0}
|\widehat Q_M(m,n)|^2
(m+1)^{2s}(n+1)^{2s}<\infty.
\]
Since $m+n+1\leq(m+1)(n+1)$ for $m,n\geq0$, our choice of $s$ gives
\[
\begin{aligned}
w_{\boldsymbol{\alpha}}(m,n)
&\leq
(m+1)^{\alpha_1^+}(n+1)^{\alpha_2^+}
(m+n+1)^{\alpha_3^+}\\
&\leq
(m+1)^{\alpha_1^++\alpha_3^+}
(n+1)^{\alpha_2^++\alpha_3^+}\\
&\leq
(m+1)^{2s}(n+1)^{2s}.
\end{aligned}
\]
Consequently,
\[
\|Q_M\|_{\boldsymbol{\alpha}}^2
=
\sum_{m,n\geq0}
w_{\boldsymbol{\alpha}}(m,n)|\widehat Q_M(m,n)|^2
<\infty,
\]
and $Q_M\in\mathscr{D}_{\boldsymbol{\alpha}}$.
\end{proof}

The following theorem characterizes cyclicity in
$\mathscr{D}_{\boldsymbol{\alpha}}$ for polynomials with
a nonempty finite zero set on $\mathbb T^2$.

\begin{theorem}\label{fzc}
Let $\boldsymbol{\alpha}\in\mathbb R^3$,
and let $A,B,S$ be as in \eqref{eq:effective-parameters}.
Let $p\in\mathbb C[z_1,z_2]$ be irreducible, depend on both
variables, and have no zeros in $\mathbb D^2$.
Suppose that
\[
0<\#\bigl(\mathcal Z(p)\cap\mathbb T^2\bigr)<\infty.
\]
Then $p$ is cyclic in $\mathscr{D}_{\boldsymbol{\alpha}}$
if and only if $A\leq1$ or $B\leq1$ or $S\leq2$.
\end{theorem}

\begin{proof}
If $A>1$, $B>1$, and $S>2$, then Corollary~\ref{ABS}
implies that $p$ is noncyclic, since it has a zero on
$\mathbb T^2$. This proves necessity.

Conversely, suppose that $A\leq1$ or $B\leq1$ or $S\leq2$.
Write
\[
\mathcal Z(p)\cap\mathbb T^2
=
\{(\zeta_j,\eta_j):j=1,\ldots,N\}
\]
and define
\[
L_j(z_1,z_2)
=
2-\overline{\zeta_j}z_1-\overline{\eta_j}z_2,
\qquad j=1,\ldots,N.
\]
By Theorem~\ref{thm:model-polynomial-cyclicity}(iv),
each $L_j$ is cyclic in $\mathscr{D}_{\boldsymbol{\alpha}}$.
Moreover, Lemma~\ref{lem:smooth-quotient} gives a positive
integer $M$ such that
\[
Q=\frac{\prod_{j=1}^{N}L_j^M}{p}
\in\mathscr{D}_{\boldsymbol{\alpha}}.
\]
Polynomial approximation of $Q$ and boundedness of
multiplication by $p$ imply that
\[
\prod_{j=1}^{N}L_j^M=pQ\in[p]_{\boldsymbol{\alpha}}.
\]
A finite product of cyclic polynomial multipliers is cyclic,
so $\prod_{j=1}^{N}L_j^M$ is cyclic.
Thus the closed invariant subspace $[p]_{\boldsymbol{\alpha}}$
contains a cyclic polynomial, and hence $p$ is cyclic.
\end{proof}

 \begin{remark}
    Note that if $\max\{A,B\} \leq 1,$ the we can use the model polynomials $z_i-\zeta_i,\, \zeta_i \in \mathbb T$ to prove the cyclicity of $p$ in Theorem~\ref{fzc}. 
\end{remark}

\section{Polynomials having infinitely many zeros on
\texorpdfstring{$\mathbb T^2$}{the two-torus}}
\label{i-b-z}

We prove noncyclicity of a polynomial $p$ in $\mathscr{D}_{\boldsymbol{\alpha}}$
when $\alpha_1+\alpha_2+\alpha_3>1$ by constructing a bounded linear functional
supported on a smooth arc of zeros of $p$ on $\mathbb T^2$.
This functional annihilates every polynomial multiple of $p$
but does not vanish on the constant function $1$.

\medskip

We first recall the geometric argument for the existence
of a suitable arc from
\cite[Section~2.2 and proof of Theorem~5.2]{BKKLSS2016}.
An essential additional step is to show that its angular parametrization
$t\mapsto(e^{it},e^{i\varphi(t)})$ can be chosen with
$-C\leq\varphi'(t)\leq-c<0$.
In particular, the derivative is uniformly bounded away
from zero.
These bounds allow us to combine integration by parts
with Parseval's identity to establish boundedness of
the functional when $\alpha_1+\alpha_2+\alpha_3>1$.

\subsection{Local geometry of the toral zero set}

The following proposition establishes a geometric property
of polynomial zero sets that will be used in the proof
of the main result of this section.

\begin{proposition}\label{prop:toral-zero-arc}
Let $p\in\mathbb C[z_1,z_2]$ be irreducible, depend on both
variables, have no zeros in $\mathbb D^2$, and have infinitely
many zeros on $\mathbb T^2$.
Then there exist a nonempty open interval $I\subset\mathbb R$,
a function $\varphi\in C^\infty(I;\mathbb R)$, and constants
$c,C>0$ such that:
\begin{enumerate}
\item
The smooth arc
\[
\Gamma=\{(e^{it},e^{i\varphi(t)}):t\in I\}
\]
is contained in $\mathcal Z(p)\cap\mathbb T^2$.

\item
The derivative of $\varphi$ satisfies
\begin{equation}\label{eq:toral-arc-slope}
-C\leq\varphi'(t)\leq-c<0,
\qquad t\in I.
\end{equation}
\end{enumerate}
\end{proposition}

\begin{proof}
Part~(i) follows from the geometric discussion in
\cite[Section~2.2 and proof of Theorem~5.2]{BKKLSS2016}.
For the reader's convenience, we include the argument,
also ensuring that the angular derivative is nonzero.

The set $\mathcal Z(p)\cap\mathbb T^2$ is semialgebraic.
By \cite[Proposition~2.9.10]{BCR1998}, it is a finite
disjoint union of connected Nash submanifolds, which are
real analytic.
Since this set is infinite, at least one of these submanifolds
has positive dimension and therefore contains a regular
real-analytic arc.
Moreover, since $p$ is irreducible and depends on both
variables, each partial derivative $p_{z_j}$ is nonzero,
has smaller degree than $p$, and is therefore relatively
prime to $p$.
By B\'ezout's theorem, the sets
$\mathcal Z(p)\cap\mathcal Z(p_{z_j})$, $j=1,2$, are finite.
We may thus restrict the arc so that both partial derivatives
are nonzero at every point of it.

Using local angular coordinates, write the arc as
\[
\gamma(s)
=\bigl(e^{i\theta_1(s)},e^{i\theta_2(s)}\bigr),
\qquad s\in J,
\]
where $J\subset\mathbb R$ is a nonempty open interval,
$\theta_1,\theta_2$ are real analytic, and
$\gamma'(s)\neq0$.
Differentiating $p(\gamma(s))=0$ gives
\[
p_{z_1}(\gamma(s))\,i e^{i\theta_1(s)}\theta_1'(s)
+
p_{z_2}(\gamma(s))\,i e^{i\theta_2(s)}\theta_2'(s)
=0.
\]
Since both partial derivatives are nonzero, the vanishing
of either $\theta_1'(s)$ or $\theta_2'(s)$ would force
the other to vanish, contradicting $\gamma'(s)\neq0$.
Thus both angular derivatives are nonzero.

Since $\theta_1'$ is continuous and nonzero on $J$,
$\theta_1$ is strictly monotone and has a real-analytic
inverse on $I_0=\theta_1(J)$.
We may therefore use $t=\theta_1(s)$ as a parameter.
Setting $\varphi=\theta_2\circ\theta_1^{-1}$, we obtain
\[
\{(e^{it},e^{i\varphi(t)}):t\in I_0\}
\subseteq\mathcal Z(p)\cap\mathbb T^2,
\]
where $\varphi$ is real analytic and
\[
\varphi'(\theta_1(s))
=\frac{\theta_2'(s)}{\theta_1'(s)}
\neq0,
\qquad s\in J.
\]
This establishes~(i).

To prove~(ii), we first show that $\varphi'<0$.
Since $p_{z_2}\neq0$ on the arc, the holomorphic implicit
function theorem, after restricting $I_0$ if necessary,
provides a nonvanishing holomorphic function $h$ on a
neighborhood of $\{e^{it}:t\in I_0\}$ such that
\[
p(z,h(z))=0,
\qquad h(e^{it})=e^{i\varphi(t)}.
\]
For each fixed $t\in I_0$ and $r<1$ sufficiently close
to $1$, we have
\[
p(re^{it},h(re^{it}))=0.
\]
Since $|re^{it}|<1$ and $p$ has no zeros in $\mathbb D^2$,
it follows that $|h(re^{it})|\geq1$.
On the other hand, $|h(e^{it})|=1$.
Consequently,
\[
\left.\frac{\partial}{\partial r}
\log|h(re^{it})|\right|_{r=1}
=
\lim_{r\uparrow1}
\frac{\log|h(re^{it})|-\log|h(e^{it})|}{r-1}
\leq0.
\]
Here the derivative exists because $h$ is holomorphic
and nonvanishing near $e^{it}$.

Differentiating the identity
$h(e^{it})=e^{i\varphi(t)}$ gives
\[
ie^{it}h'(e^{it})
=i\varphi'(t)h(e^{it}).
\]
Dividing by $ih(e^{it})$ and taking real parts, we obtain
\[
\varphi'(t)
=\operatorname{Re}\frac{e^{it}h'(e^{it})}{h(e^{it})}
=\left.\frac{\partial}{\partial r}
  \log|h(re^{it})|\right|_{r=1}
\leq0.
\]
Since $\varphi'(t)\neq0$, we conclude that
$\varphi'(t)<0$ throughout $I_0$.

Finally, choose a nonempty open interval $I$ whose closure
is a compact subset of $I_0$.
Continuity and strict negativity of $\varphi'$ give
\[
c=\min_{t\in\overline I}\bigl(-\varphi'(t)\bigr)>0,
\qquad
C=\max_{t\in\overline I}\bigl(-\varphi'(t)\bigr)<\infty.
\]
Restricting $\varphi$ to $I$, these constants satisfy
\eqref{eq:toral-arc-slope}, proving~(ii).
\end{proof}

\begin{remark}
Let $I,\varphi,c,C$ be as in
Proposition~\ref{prop:toral-zero-arc}.
By part~(ii), for $m,n\geq0$, the derivative of
$t\mapsto mt+n\varphi(t)$ can vanish at a point of $I$
only if
\[
cn\leq m\leq Cn.
\]
In this region,
$m+1\asymp n+1\asymp m+n+1$, and consequently
\[
w_{\boldsymbol{\alpha}}(m,n)\asymp(n+1)^{\alpha_1+\alpha_2+\alpha_3}.
\]
This comparison will be used in the proof of
Lemma~\ref{lem:curve-functional} below to establish boundedness
of the functional supported on the arc when $\alpha_1+\alpha_2+\alpha_3>1$.
\end{remark}

\subsection{A bounded functional supported on an arc}

\begin{lemma}\label{lem:off-cone-decay}
Let $I\subset\mathbb R$ be an open interval, and let
$\varphi\in C^\infty(I;\mathbb R)$ satisfy
$-C\leq\varphi'\leq-c<0$ for some $c,C>0$.
For $\chi\in C_c^\infty(I)$, set
\[
c_{m,n}=\int_I e^{i(mt+n\varphi(t))}\chi(t)\,dt,
\qquad m,n\geq0.
\]
Choose $0<c_0<c$ and $C_0>C$, and define
\[
\mathcal C
=\{(m,n)\in\mathbb N_0^2:c_0n\leq m\leq C_0n\}.
\]
For every integer $N\geq1$, there is a constant $K_N>0$ such that
\[
|c_{m,n}|\leq K_N(1+m+n)^{-N},
\qquad (m,n)\notin\mathcal C.
\]
\end{lemma}

\begin{proof}
We divide the argument into three claims.
The first gives a lower bound for the derivative outside
the cone $\mathcal C$.

\medskip
\noindent\textbf{Claim 1.}
For every $(m,n)\in\mathbb N_0^2\setminus\mathcal C$ and $t\in I$,
\[
|m+n\varphi'(t)|\geq\delta(m+n),
\]
where
\[
\delta=\min\left\{
\frac{c-c_0}{1+c_0},
\frac{C_0-C}{1+C_0}
\right\}>0
\]
is independent of $m,n$, and $t$.

\smallskip

\noindent \textit{Proof of Claim 1.}
If $(m,n)\notin\mathcal C$, then either $m<c_0n$ or $m>C_0n$.

If $m<c_0n$, then $m+n\varphi'(t)<0$.
Since $m+n\leq(1+c_0)n$, we obtain
\[
\begin{aligned}
|m+n\varphi'(t)|
&=-m-n\varphi'(t)\\
&\geq cn-m\\
&\geq(c-c_0)n\\
&\geq\frac{c-c_0}{1+c_0}(m+n).
\end{aligned}
\]
If $m>C_0n$, then $m+n\varphi'(t)>0$.
Since $n\leq(m+n)/(1+C_0)$, we obtain
\[
\begin{aligned}
|m+n\varphi'(t)|
&=m+n\varphi'(t)\\
&\geq m-Cn\\
&=(m+n)-(C+1)n\\
&\geq\frac{C_0-C}{1+C_0}(m+n).
\end{aligned}
\]
The claim follows from the definition of $\delta$.\hfill$\blacksquare$

\medskip

For $(m,n)\in\mathbb N_0^2\setminus\mathcal C$, set
\[
\rho=m+n,\qquad
\Psi_{m,n}(t)=\frac{mt+n\varphi(t)}{\rho}.
\]
Since $(0,0)\in\mathcal C$, we have $\rho\geq1$.
We next establish derivative bounds for $\Psi_{m,n}$
that are uniform in $m,n$.

\medskip
\noindent\textbf{Claim 2.}
For every $(m,n)\in\mathbb N_0^2\setminus\mathcal C$, we have
\[
|\Psi_{m,n}'(t)|\geq\delta,
\qquad t\in I.
\]
Moreover, there are constants $B_k,D_k>0$, independent of
$m,n$, such that
\[
\begin{aligned}
\sup_{t\in\operatorname{supp}\chi}|\Psi_{m,n}^{(k)}(t)|
&\leq B_k,
&& k\geq1,\\
\sup_{t\in\operatorname{supp}\chi}
\left|
\frac{d^k}{dt^k}\left(\frac1{\Psi_{m,n}'(t)}\right)
\right|
&\leq D_k,
&& k\geq0.
\end{aligned}
\]

\smallskip
\noindent\textit{Proof of Claim 2.}
By Claim~1,
\[
|\Psi_{m,n}'(t)|
=\frac{|m+n\varphi'(t)|}{m+n}
\geq\delta,
\qquad t\in I.
\]
Differentiating $\Psi_{m,n}$ gives
\[
\begin{aligned}
\Psi_{m,n}'(t)
&=\frac{m}{m+n}+\frac{n}{m+n}\varphi'(t),\\
\Psi_{m,n}^{(k)}(t)
&=\frac{n}{m+n}\varphi^{(k)}(t),
\qquad k\geq2.
\end{aligned}
\]
Since $0\leq m/(m+n),n/(m+n)\leq1$, we have
$|\Psi_{m,n}'(t)|\leq1+C$ and
$|\Psi_{m,n}^{(k)}(t)|\leq|\varphi^{(k)}(t)|$ for $k\geq2$.
Each derivative of $\varphi$ is continuous and therefore
bounded on the compact set $\operatorname{supp}\chi$.
Thus, for every integer $k\geq1$, there is a constant
$B_k>0$, independent of $m,n$, such that
\[
\sup_{t\in\operatorname{supp}\chi}
|\Psi_{m,n}^{(k)}(t)|
\leq B_k.
\]

The bound $|\Psi_{m,n}'|\geq\delta$ gives
\[
\sup_{t\in\operatorname{supp}\chi}
\left|\frac1{\Psi_{m,n}'(t)}\right|
\leq\frac1\delta=:D_0.
\]
For every integer $k\geq1$, repeated differentiation gives
\[
\frac{d^k}{dt^k}\left(\frac1{\Psi_{m,n}'}\right)
=
\sum_{\ell=1}^{k}
\sum_{\substack{j_1+\cdots+j_\ell=k\\j_1,\ldots,j_\ell\geq1}}
\frac{(-1)^\ell k!}{j_1!\cdots j_\ell!}
\frac{\prod_{r=1}^{\ell}\Psi_{m,n}^{(j_r+1)}}
     {(\Psi_{m,n}')^{\ell+1}}.
\]
On $\operatorname{supp}\chi$, each numerator is bounded by
$\prod_{r=1}^{\ell}B_{j_r+1}$, while the absolute value
of the denominator is at least $\delta^{\ell+1}$.
Since the sum is finite for each fixed $k$, there is a constant
$D_k>0$, independent of $m,n$, such that
\[
\sup_{t\in\operatorname{supp}\chi}
\left|
\frac{d^k}{dt^k}\left(\frac1{\Psi_{m,n}'(t)}\right)
\right|
\leq D_k.
\]
Together with the estimate for $D_0$, this proves Claim~2.
\hfill$\blacksquare$

\medskip

The final claim gives an integral representation obtained
by repeated integration by parts.

\medskip
\noindent\textbf{Claim 3.}
For every integer $N\geq1$ and every
$(m,n)\in\mathbb N_0^2\setminus\mathcal C$, there is a function
$b_N\in C_c^\infty(I)$, supported in $\operatorname{supp}\chi$,
such that
\[
c_{m,n}
=\rho^{-N}\int_I e^{i\rho\Psi_{m,n}(t)}b_N(t)\,dt,
\qquad
\int_I|b_N(t)|\,dt\leq C_N.
\]
Here $b_N$ may depend on $m,n$, but the constant $C_N$ does not.

\smallskip
\noindent\textit{Proof of Claim 3.}
Fix $(m,n)\in\mathbb N_0^2\setminus\mathcal C$.
Using
\[
e^{i\rho\Psi_{m,n}(t)}
=\frac{1}{i\rho\Psi_{m,n}'(t)}
\frac{d}{dt}e^{i\rho\Psi_{m,n}(t)},
\]
integration by parts gives
\[
\begin{aligned}
c_{m,n}
&=\frac1\rho\int_I
\frac{\chi(t)}{i\Psi_{m,n}'(t)}
\frac{d}{dt}e^{i\rho\Psi_{m,n}(t)}\,dt\\
&=-\frac1\rho\int_I e^{i\rho\Psi_{m,n}(t)}
\frac{d}{dt}\left(\frac{\chi(t)}{i\Psi_{m,n}'(t)}\right)\,dt.
\end{aligned}
\]
There is no boundary term because $\chi$ is compactly
supported in $I$.

Define recursively
\[
b_0=\chi,
\qquad
b_{j+1}
=-\frac{d}{dt}\left(\frac{b_j}{i\Psi_{m,n}'}\right),
\qquad j\geq0.
\]
Each $b_j$ is smooth and supported in $\operatorname{supp}\chi$,
so no boundary terms arise in subsequent integrations.
After $N$ integrations by parts, we obtain
\[
c_{m,n}
=\rho^{-N}\int_I e^{i\rho\Psi_{m,n}(t)}b_N(t)\,dt.
\]

It remains to bound the integral of $|b_N|$ uniformly in $m,n$.
Set $K=\operatorname{supp}\chi$.
By the Leibniz rule and the bounds $D_\ell$ from Claim~2,
\[
\sup_{t\in K}|b_{j+1}^{(k)}(t)|
\leq
\sum_{\ell=0}^{k+1}\binom{k+1}{\ell}
D_\ell
\sup_{t\in K}|b_j^{(k+1-\ell)}(t)|,
\qquad j,k\geq0.
\]
Since $b_0=\chi$ is fixed and all its derivatives are bounded
on $K$, induction on $j$ shows that every derivative of
$b_j$ is uniformly bounded on $K$, independently of $m,n$.
In particular, there is a constant $E_N>0$ such that
\[
\sup_{t\in K}|b_N(t)|\leq E_N.
\]
Because $b_N$ is supported in $K$, we conclude that
\[
\int_I|b_N(t)|\,dt
\leq |K|E_N=:C_N,
\]
where $|K|$ denotes the Lebesgue measure of $K$.
The constant $C_N$ is independent of $m,n$.
\hfill$\blacksquare$

\medskip

Finally, since $\rho=m+n\geq1$ outside $\mathcal C$, we have
$1+\rho\leq2\rho$.
Thus Claim~3 and $|e^{i\rho\Psi_{m,n}(t)}|=1$ give
\[
|c_{m,n}|
\leq C_N\rho^{-N}
\leq 2^NC_N(1+\rho)^{-N}
=K_N(1+m+n)^{-N},
\]
where $K_N=2^NC_N$ is independent of $m,n$.
This completes the proof.
\end{proof}

We write $I\Subset J$ when the closure of $I$ is a compact subset
of $J$.

The following lemma shows that integration against a smooth
density with compact support on an arc satisfying
\eqref{eq:toral-arc-slope} defines a bounded linear functional
on $\mathscr{D}_{\boldsymbol{\alpha}}$ whenever $\alpha_1+\alpha_2+\alpha_3>1$.

\begin{lemma}\label{lem:curve-functional}
Let $\boldsymbol{\alpha}\in\mathbb R^3$ and suppose that
$S=\alpha_1+\alpha_2+\alpha_3>1$.
Let $I\Subset J=(a,a+2\pi)$, and let
$\varphi\in C^\infty(I;\mathbb R)$ satisfy
$-C\leq\varphi'\leq-c<0$ for some $c,C>0$.
For every $\chi\in C_c^\infty(I)$, the functional
\[
\Lambda_\chi(f)
=\int_I f(e^{it},e^{i\varphi(t)})\chi(t)\,dt,
\]
initially defined on polynomials, extends continuously to
$\mathscr{D}_{\boldsymbol{\alpha}}$.
\end{lemma}

\begin{proof}
For $m\in\mathbb Z$ and $n\geq0$, put
\[
c_{m,n}=\int_I e^{i(mt+n\varphi(t))}\chi(t)\,dt.
\]
We prove that
\begin{equation}\label{eq:cone-functional-energy}
\sum_{m,n\geq0}
\frac{|c_{m,n}|^2}{w_{\boldsymbol{\alpha}}(m,n)}<\infty.
\end{equation}
Choose $\mathcal C$ as in Lemma~\ref{lem:off-cone-decay}.
Outside $\mathcal C$, that lemma applies with any $N$.
Since
\[
w_{\boldsymbol{\alpha}}(m,n)^{-1}
\leq(1+m+n)^R,
\qquad R=\sum_{j=1}^3\max\{-\alpha_j,0\},
\]
choosing $2N-R>2$ gives
\[
\begin{aligned}
\sum_{(m,n)\in\mathbb N_0^2\setminus\mathcal C}
\frac{|c_{m,n}|^2}{w_{\boldsymbol{\alpha}}(m,n)}
&\lesssim\sum_{m,n\geq0}(1+m+n)^{R-2N}\\
&=\sum_{k\geq0}(k+1)^{R-2N+1}<\infty.
\end{aligned}
\]

Inside $\mathcal C$, we have
$m+1\asymp n+1\asymp m+n+1$, so
$$w_{\boldsymbol{\alpha}}(m,n)\asymp(n+1)^{\alpha_1+\alpha_2+\alpha_3}.$$
For fixed $n$, extend $g_n(t)=\chi(t)e^{in\varphi(t)}$
by zero from $I$ to $J$.
With the convention
$\widehat g_n(k)=(2\pi)^{-1}\int_Jg_n(t)e^{-ikt}\,dt$,
we have $c_{m,n}=2\pi\widehat g_n(-m)$.
Parseval's identity therefore gives
\[
\sum_{m\in\mathbb Z}|c_{m,n}|^2
=2\pi\int_J|g_n(t)|^2\,dt
=2\pi\|\chi\|_{L^2(I)}^2.
\]
Consequently,
\[
\begin{aligned}
\sum_{(m,n)\in\mathcal C}
\frac{|c_{m,n}|^2}{w_{\boldsymbol{\alpha}}(m,n)}
&\lesssim\sum_{n\geq0}(n+1)^{-\alpha_1-\alpha_2-\alpha_3}
\sum_{m\in\mathbb Z}|c_{m,n}|^2\\
&=2\pi\|\chi\|_{L^2(I)}^2
\sum_{n\geq0}(n+1)^{-\alpha_1-\alpha_2-\alpha_3}<\infty,
\end{aligned}
\]
because $\alpha_1+\alpha_2+\alpha_3>1$. This proves \eqref{eq:cone-functional-energy}.

For a polynomial $f=\sum a_{m,n}z_1^mz_2^n$, Cauchy--Schwarz now gives
\[
|\Lambda_\chi(f)|
=\left|\sum_{m,n\geq0}a_{m,n}c_{m,n}\right|
\leq\|f\|_{\boldsymbol{\alpha}}
\left(\sum_{m,n\geq0}
\frac{|c_{m,n}|^2}{w_{\boldsymbol{\alpha}}(m,n)}\right)^{1/2}.
\]
By density of polynomials, $\Lambda_\chi$ extends uniquely
to a bounded linear functional on
$\mathscr{D}_{\boldsymbol{\alpha}}$.
\end{proof}

\subsection{Noncyclicity}

\begin{theorem}\label{infzero}
Let $\boldsymbol{\alpha}\in\mathbb R^3$,
and let $A,B,S$ be as in \eqref{eq:effective-parameters}.
Let $p\in\mathbb C[z_1,z_2]$ be irreducible, depend on both
variables, and have no zeros in $\mathbb D^2$.
Suppose that $\mathcal Z(p)\cap\mathbb T^2$ is infinite.
If $S>1$, then $p$ is noncyclic in
$\mathscr{D}_{\boldsymbol{\alpha}}$.
\end{theorem}

\begin{proof}
By Proposition~\ref{prop:toral-zero-arc}, the set of zeros
of $p$ on $\mathbb T^2$ contains a smooth arc
\[
\Gamma=\{(e^{it},e^{i\varphi(t)}):t\in I\}
\]
satisfying \eqref{eq:toral-arc-slope}.
After restricting $I$ if necessary, we may assume that
$I\Subset(a,a+2\pi)$ for some $a\in\mathbb R$.

Choose a nonnegative function $\chi\in C_c^\infty(I)$
such that $\int_I\chi(t)\,dt=1$.
Since $S>1$, Lemma~\ref{lem:curve-functional} shows that
the functional
\[
\Lambda_\chi(f)
=\int_I f(e^{it},e^{i\varphi(t)})\chi(t)\,dt,
\]
initially defined on polynomials, extends uniquely to
a bounded linear functional on
$\mathscr{D}_{\boldsymbol{\alpha}}$.

Since $p$ vanishes on $\Gamma$, we have
\[
\Lambda_\chi(pq)=0
\quad\text{for every }q\in\mathbb C[z_1,z_2],
\qquad
\Lambda_\chi(1)=1.
\]
By continuity, $\Lambda_\chi$ vanishes on
$[p]_{\boldsymbol{\alpha}}$.
Hence $1\notin[p]_{\boldsymbol{\alpha}}$, so $p$ is
noncyclic in $\mathscr{D}_{\boldsymbol{\alpha}}$.
\end{proof}

\section{Completion of the proof of Theorem~\ref{thm:main}}
\label{sec:proof-main}

We combine the results of the preceding sections to complete
the proof of Theorem~\ref{thm:main}.

\begin{proof}
If $\mathcal Z(p)\cap\mathbb T^2=\varnothing$, then
Lemma~\ref{lem:no-side-boundary-zeros} shows that $p$ has no
zeros on $\overline{\mathbb D}^{\,2}$.
Hence $1/p$ is holomorphic on a neighborhood of the closed
bidisk and belongs to $\mathscr{D}_{\boldsymbol{\alpha}}$.
It follows that $1=p(1/p)\in[p]_{\boldsymbol{\alpha}}$,
so $p$ is cyclic for every $\boldsymbol{\alpha}\in\mathbb R^3$.

Part~(i) is Theorem~\ref{s1}.
For part~(ii), Theorem~\ref{fzc} gives cyclicity for a nonempty
finite toral zero set, and Theorem~\ref{infzero} excludes an
infinite toral zero set. The empty case was settled above.
Finally, under the hypotheses of part~(iii), we have
$A>1$, $B>1$, and $S>2$, so Corollary~\ref{ABS} excludes
every nonempty toral zero set. The converse again follows
from the empty case.
\end{proof}

\section{Final Remarks}
\label{sec:final-remarks}

We conclude with two further consequences of norm comparisons
with anisotropic Dirichlet-type spaces.
First, these comparisons give a shorter proof of
Theorem~\ref{infzero} in two parameter ranges.
We then establish a Hilbert-sum representation of
$\mathscr{D}_{\boldsymbol{\alpha}}$ when $\alpha_3\leq0$.

\begin{proposition}
Let $\boldsymbol{\alpha}=(\alpha_1,\alpha_2,\alpha_3)\in\mathbb R^3$, and set
$S=\alpha_1+\alpha_2+\alpha_3$.
Let $p\in\mathbb C[z_1,z_2]$ be irreducible, depend on both
variables, have no zeros in $\mathbb D^2$, and have infinitely
many zeros on $\mathbb T^2$.
Then the following statements hold.
\begin{enumerate}
\item
If $\alpha_3\geq0$ and $S>1$, then $p$ is noncyclic in
$\mathscr{D}_{\boldsymbol{\alpha}}$.

\item
If $\alpha_3<0$ and $S+\alpha_3>1$, then $p$ is noncyclic in
$\mathscr{D}_{\boldsymbol{\alpha}}$.
\end{enumerate}
\end{proposition}

\begin{proof}
Suppose first that $\alpha_3\geq0$. Fix
$t\in[0,\alpha_3]$ and set
\[
\beta_1=\alpha_1+t,
\qquad
\beta_2=\alpha_2+\alpha_3-t.
\]
Since
\[
(m+n+1)^{\alpha_3}
\geq
(m+1)^t(n+1)^{\alpha_3-t},
\]
we obtain
\[
w_{\boldsymbol{\alpha}}(m,n)
\geq
(m+1)^{\beta_1}(n+1)^{\beta_2}.
\]
Thus
\[
\mathscr{D}_{\boldsymbol{\alpha}}
\hookrightarrow
\mathfrak D_{(\beta_1,\beta_2)}.
\]
Since $\beta_1+\beta_2=S>1$ and $p$ has infinitely many
zeros on $\mathbb T^2$,
Theorem~\ref{thm:anisotropic-cyclicity}(ii)--(iii)
implies that $p$ is noncyclic in
$\mathfrak D_{(\beta_1,\beta_2)}$.
If $p$ were cyclic in $\mathscr{D}_{\boldsymbol{\alpha}}$,
the continuous inclusion would make it cyclic in
$\mathfrak D_{(\beta_1,\beta_2)}$, a contradiction.

Suppose next that $\alpha_3<0$. Since
\[
m+n+1\leq(m+1)(n+1),
\]
we have
\[
(m+n+1)^{\alpha_3}
\geq
(m+1)^{\alpha_3}(n+1)^{\alpha_3}.
\]
Consequently,
\[
\mathscr{D}_{\boldsymbol{\alpha}}
\hookrightarrow
\mathfrak D_{(\alpha_1+\alpha_3,\alpha_2+\alpha_3)}.
\]
The sum of the two anisotropic parameters is
\[
\alpha_1+\alpha_2+2\alpha_3=S+\alpha_3.
\]
If $S+\alpha_3>1$, then
Theorem~\ref{thm:anisotropic-cyclicity}(ii)--(iii)
implies that $p$ is noncyclic in
$\mathfrak D_{(\alpha_1+\alpha_3,\alpha_2+\alpha_3)}$.
The continuous inclusion again implies that $p$ is
noncyclic in $\mathscr{D}_{\boldsymbol{\alpha}}$.
\end{proof}

For $\alpha_3\leq0$, the continuous inclusions in
Lemma~\ref{lem:mixed-anisotropic-comparison}(i)
can be strengthened to the following Hilbert-sum
representation.

\begin{proposition}\label{prop:hilbert-sum}
Let $\boldsymbol{\alpha}=(\alpha_1,\alpha_2,\alpha_3)\in\mathbb R^3$, and set
$A=\alpha_1+\alpha_3$ and $B=\alpha_2+\alpha_3$.
If $\alpha_3\leq0$, then
\[
\mathscr{D}_{\boldsymbol{\alpha}}
=
\mathfrak D_{(A,\alpha_2)}
+
\mathfrak D_{(\alpha_1,B)}
\]
with equivalent norms:
\[
\|f\|_{\boldsymbol{\alpha}}^2
\asymp
\inf_{\substack{f=g+h\\
g\in\mathfrak D_{(A,\alpha_2)}\\
h\in\mathfrak D_{(\alpha_1,B)}}}
\left(
\|g\|_{\mathfrak D_{(A,\alpha_2)}}^2
+
\|h\|_{\mathfrak D_{(\alpha_1,B)}}^2
\right).
\]
\end{proposition}

\begin{proof}
Put
\[
u_{m,n}=(m+1)^A(n+1)^{\alpha_2},
\qquad
v_{m,n}=(m+1)^{\alpha_1}(n+1)^B,
\]
and let $r=-\alpha_3\geq0$.
The inequalities
\[
\max\{m+1,n+1\}
\leq m+n+1
\leq2\max\{m+1,n+1\}
\]
imply that
\[
(m+n+1)^r
\asymp
(m+1)^r+(n+1)^r,
\]
with constants independent of $m,n$.
This comparison also holds when $r=0$.
Multiplying by $(m+1)^{-\alpha_1}(n+1)^{-\alpha_2}$,
we obtain
\[
\begin{aligned}
w_{\boldsymbol{\alpha}}(m,n)^{-1}
&=(m+1)^{-\alpha_1}(n+1)^{-\alpha_2}
  (m+n+1)^r\\
&\asymp
(m+1)^{-\alpha_1+r}(n+1)^{-\alpha_2}
+
(m+1)^{-\alpha_1}(n+1)^{-\alpha_2+r}\\
&=(m+1)^{-A}(n+1)^{-\alpha_2}
+
(m+1)^{-\alpha_1}(n+1)^{-B}\\
&=u_{m,n}^{-1}+v_{m,n}^{-1}.
\end{aligned}
\]
Here we used $A=\alpha_1-r$ and $B=\alpha_2-r$.
Therefore
\[
w_{\boldsymbol{\alpha}}(m,n)
\asymp
\frac{u_{m,n}v_{m,n}}{u_{m,n}+v_{m,n}}.
\]

We next establish the scalar minimization identity underlying
the Hilbert-sum norm.
For $u,v>0$ and $a,b,c\in\mathbb C$ with $a+b=c$,
completing the square gives
\[
u|a|^2+v|b|^2
=
(u+v)\left|a-\frac{v}{u+v}c\right|^2
+\frac{uv}{u+v}|c|^2.
\]
The first term is nonnegative and vanishes precisely when
$a=vc/(u+v)$.
Consequently,
\[
\min_{a+b=c}\bigl(u|a|^2+v|b|^2\bigr)
=
\frac{uv}{u+v}|c|^2,
\]
with the minimum attained at
\[
a=\frac{v}{u+v}c,
\qquad
b=\frac{u}{u+v}c.
\]

Suppose now that $f=g+h$, where
\[
\begin{aligned}
g&=\sum_{m,n\geq0}a_{m,n}z_1^mz_2^n
\in\mathfrak D_{(A,\alpha_2)},\\
h&=\sum_{m,n\geq0}b_{m,n}z_1^mz_2^n
\in\mathfrak D_{(\alpha_1,B)}.
\end{aligned}
\]
The Taylor coefficients of $f$ are
$c_{m,n}=a_{m,n}+b_{m,n}$.
Applying the scalar inequality to each coefficient and
summing, we find
\[
\begin{aligned}
\sum_{m,n\geq0}
\frac{u_{m,n}v_{m,n}}{u_{m,n}+v_{m,n}}
|c_{m,n}|^2
&\leq
\sum_{m,n\geq0}
\left(
u_{m,n}|a_{m,n}|^2+v_{m,n}|b_{m,n}|^2
\right)\\
&=
\|g\|_{\mathfrak D_{(A,\alpha_2)}}^2
+
\|h\|_{\mathfrak D_{(\alpha_1,B)}}^2
<\infty.
\end{aligned}
\]
By the weight comparison, $f\in\mathscr{D}_{\boldsymbol{\alpha}}$.
This proves that the sum of the two anisotropic spaces
is contained in $\mathscr{D}_{\boldsymbol{\alpha}}$.

Conversely, let
\[
f=\sum_{m,n\geq0}c_{m,n}z_1^mz_2^n
\in\mathscr{D}_{\boldsymbol{\alpha}},
\]
and choose the coefficientwise minimizers
\[
a_{m,n}
=\frac{v_{m,n}}{u_{m,n}+v_{m,n}}c_{m,n},
\qquad
b_{m,n}
=\frac{u_{m,n}}{u_{m,n}+v_{m,n}}c_{m,n}.
\]
Since $u_{m,n},v_{m,n}>0$, we have
$|a_{m,n}|\leq|c_{m,n}|$ and
$|b_{m,n}|\leq|c_{m,n}|$.
Thus the power series
\[
g=\sum_{m,n\geq0}a_{m,n}z_1^mz_2^n,
\qquad
h=\sum_{m,n\geq0}b_{m,n}z_1^mz_2^n
\]
converge absolutely and locally uniformly on $\mathbb D^2$,
as does the Taylor series of $f$.
They therefore define holomorphic functions, and
$a_{m,n}+b_{m,n}=c_{m,n}$ gives $f=g+h$.

For this choice of coefficients, the scalar minimum is
attained at every index, so
\[
\begin{aligned}
\sum_{m,n\geq0}
\left(
u_{m,n}|a_{m,n}|^2+v_{m,n}|b_{m,n}|^2
\right)
&=
\sum_{m,n\geq0}
\frac{u_{m,n}v_{m,n}}{u_{m,n}+v_{m,n}}
|c_{m,n}|^2\\
&\asymp
\|f\|_{\boldsymbol{\alpha}}^2
<\infty.
\end{aligned}
\]
Hence $g\in\mathfrak D_{(A,\alpha_2)}$ and
$h\in\mathfrak D_{(\alpha_1,B)}$, proving the reverse
inclusion.

Finally, the lower bound for arbitrary decompositions
and the explicit minimizing decomposition give
\[
\begin{aligned}
&\inf_{\substack{f=g+h\\
g\in\mathfrak D_{(A,\alpha_2)}\\
h\in\mathfrak D_{(\alpha_1,B)}}}
\left(
\|g\|_{\mathfrak D_{(A,\alpha_2)}}^2
+
\|h\|_{\mathfrak D_{(\alpha_1,B)}}^2
\right)\\
&\qquad=
\sum_{m,n\geq0}
\frac{u_{m,n}v_{m,n}}{u_{m,n}+v_{m,n}}
|c_{m,n}|^2
\asymp
\|f\|_{\boldsymbol{\alpha}}^2.
\end{aligned}
\]
This proves the asserted equality of spaces and equivalence
of norms.
\end{proof}

\subsection*{Data Availability} No data were used or generated in this study.

\subsection*{AI Disclosure}
During the preparation of this work, the authors used OpenAI's GPT-6 to explore possible proof strategies
and organize intermediate \LaTeX{} drafts. The authors subsequently reviewed and verified all mathematical arguments and take full responsibility for the content of the manuscript.

\end{document}